\documentclass[12pt]{article}

\newif\iffrench\frenchfalse
\usepackage[utf8]{inputenc}
\usepackage[T1]{fontenc}

\iffrench
\usepackage[french]{babel}
\else
\usepackage[english]{babel}
\fi

\usepackage{pslatex}
\usepackage{lmodern}

\usepackage{amssymb}
\usepackage{amsthm}
\usepackage{bbold} % indicatrices
\usepackage{dsfont}
\usepackage{mathtools}
\usepackage{thmtools}

\usepackage{varioref} % Intelligent page references
\usepackage{hyperref} % Extensive support for hypertext in LaTeX
\hypersetup{
    colorlinks=true,
    linkcolor=blue,
    }

\usepackage{tikz}
\usepackage{pstricks-add}

\usepackage{enumitem} % Control layout of itemize, enumerate, description
\usepackage[babel=true]{csquotes} % Context sensitive quotation facilities
\usepackage{color}
\usepackage{setspace}
\usepackage{array}
\usepackage{layout}

\usepackage{stmaryrd} % St Mary Road symbols for theoretical computer science
\usepackage{mathrsfs} % jolies cursives (mathscr)

\usepackage[margin=20pt,font=small,labelfont=bf]{caption} % jolis titres de figures

\usepackage[
  backend=biber,
  citestyle=numeric-comp,
  giveninits=true,
  url=false,
  isbn=false,
  doi=false,
  eprint=false,
  maxbibnames=50,
]{biblatex}
\newcommand{\smarge}[2]{\usepackage[top=#1,bottom=#1+1cm,left=#1-#2,right=#1]{geometry}}

\iffrench

\addto\captionsfrench{} % nom des tableaux en français

\newcommand{\ioo}[1]{\left]#1\right[}
\newcommand{\ifo}[1]{\left[#1\right[}

\newtheorem{thm}{Theorème}[section]
\newtheorem{ppn}[thm]{Proposition}
\newtheorem{cor}[thm]{Corollaire}
\newtheorem{lem}[thm]{Lemme}
\newtheorem{dfi}[thm]{Définition}
\newtheorem{cjt}[thm]{Conjecture}
\newcommand{\pth}[1]{\left(#1\right)}
\newcommand{\cro}[1]{\left[#1\right]}
\newcommand{\acc}[1]{\left\{#1\right\}}
\newcommand{\abs}[1]{\left|#1\right|}
\newcommand{\dabs}[1]{\left\|#1\right\|}

\newcommand{\floor}[1]{\left\lfloor#1\right\rfloor}
\newcommand{\ceil}[1]{\left\lceil#1\right\rceil}

\else

\newtheorem{thm}{Theorem}[section]
\newtheorem{ppn}[thm]{Proposition}

\newtheorem{lem}[thm]{Lemma}
\newtheorem{dfi}[thm]{Definition}

\theoremstyle{remark}
\newtheorem{req}[thm]{Remark}

\newcommand{\ioo}[1]{\left(#1\right)}
\newcommand{\ifo}[1]{\left[#1\right)}

\newcommand{\pth}[1]{\left(#1\right)}
\newcommand{\cro}[1]{\left[#1\right]}
\newcommand{\acc}[1]{\left\{#1\right\}}
\newcommand{\abs}[1]{\left|#1\right|}
\newcommand{\dabs}[1]{\left\|#1\right\|}

\newcommand{\floor}[1]{\left\lfloor#1\right\rfloor}
\newcommand{\ceil}[1]{\left\lceil#1\right\rceil}

\newcommand{\ie}{i.e. } 
\newcommand{\eg}{e.g. }

\fi
\newcommand{\esp}{\hspace{1cm}}
\newcommand{\comment}[1]{\hspace{0.5cm}\text{#1}\hspace{0cm}}

\newcommand{\tq}{\hspace{0.25cm}/ \hspace{0.25cm}}
\newcommand{\vg}{,\,}

\newcommand{\goq}{\geqslant}
\newcommand{\loq}{\leqslant}
\newcommand{\eps}{\varepsilon}
\newcommand{\ind}{\mathbb{1}}

\newcommand{\fig}[3]{\begin{figure}[ht]\begin{center}\includegraphics[width=#1cm]{#2}\end{center}\caption{#3}\end{figure}}

\newcommand{\sys}[1]{\begin{equation}\left\{\begin{aligned}#1\end{aligned}\right.\end{equation}}
\newcommand{\syss}[1]{\begin{equation*}\left\{\begin{aligned}#1\end{aligned}\right.\end{equation*}}

\newcommand{\pmat}[1]{\begin{pmatrix} #1\end{pmatrix}}

\newcommand{\Tr}{\text{Tr}}

\newcommand{\de}{\,\mathrm{d}}

\newcommand{\dr}{\partial}

\newcommand{\Er}{\mathds{R}}
\newcommand{\Zed}{\mathds{Z}}

\newcommand{\Te}{\mathds{T}}

\newcommand{\prb}[1]{\mathds{P}\pth{#1}}
\newcommand{\psf}{\mathds{P}} % proba starting from
\newcommand{\Esp}[1]{\mathds{E}\cro{#1}}
\newcommand{\esf}{\mathds{E}} % expectation starting from

\newcommand{\kt }{\,|\,}

\newcommand{\mca}{\mathcal{A}}
\newcommand{\mcb}{\mathcal{B}}
\newcommand{\mcc}{\mathcal{C}}
\newcommand{\mcd}{\mathcal{D}}
\newcommand{\mce}{\mathcal{E}}

\newcommand{\mcg}{\mathcal{G}}

\newcommand{\mcl}{\mathcal{L}}
\newcommand{\mcm}{\mathcal{M}}

\newcommand{\mcw}{\mathcal{W}}

\newcommand{\mcy}{\mathcal{Y}}
\newcommand{\mcz}{\mathcal{Z}}

\smarge{2.5cm}{0cm}

\renewbibmacro{in:}{}

\ifpdf
\hypersetup{
  pdftitle={The principal eigenvalue problem for a strongly anisotropic second-order elliptic operator}
  pdfauthor={N. Boutillon}
}
\fi

\usepackage{soul}

\setstcolor{violet}

\title{The principal eigenvalue problem for a strongly anisotropic second-order elliptic operator
}

\author{Nathanaël Boutillon \\
\footnotesize{INRAE, BioSP, 84914, Avignon, France}\\
\footnotesize{Aix Marseille Univ, CNRS, I2M, Marseille, France}\\
\footnotesize{{nathanael.boutillon@inrae.fr}}\\
}

\date{}

\begin{document}

%\linenumbers

\maketitle

\begin{abstract}
  We consider an elliptic operator in which the second-order term is very small in one direction. In this regime, we study the behaviour of the principal eigenfunction and of the principal eigenvalue.
  Our first result deals with the limit of the principal eigenfunction and is shown with a representation of the principal eigenfunction as a quasi-stationary distribution. Subsequent results deal with the limit of the principal eigenvalue and are shown using Hamilton-Jacobi equations.
\end{abstract}



\tableofcontents

%\newpage
\vspace{1cm}

\section{Introduction}\label{s:intro}

\subsection{Presentation of the problem}
Let $n\goq 1$ and let $\mcy$ be either the $n$-dimensional torus $\Te^n=\Er^n/\Zed^n$, or an open bounded domain in $\Er^{n}$ with $\mcc^{2,\alpha}$ boundary, for some fixed $\alpha\in\ioo{0,1}$. Likewise, let $p\goq 1$ and let $\mcz$ be either the $p$-dimensional torus $\Te^p=\Er^p/\Zed^p$, or an open bounded domain in~$\Er^{p}$ with $\mcc^{2,\alpha}$ boundary. We assume that at least~$\mcy$ or~$\mcz$ is a torus. 
Let $\mcw:=\mcy\times\mcz$. If~$\mcy$ or~$\mcz$ is not a torus, the boundary of $\mcw$ is nonempty, has regularity $\mcc^{2,\alpha}$, and we let $\nu$ be the outward unit normal on $\dr\mcw$.

Let us focus on the principal eigenvalue problem
\sys{\label{eq:intro_base}
  \mcl_y\varphi+\mcl_z\varphi+c(y,z)\varphi&=k\varphi&\text{in $\mcw$},\\
  \nu\cdot\nabla\varphi&=0&\text{on $\dr\mcw$},\\
  \varphi&>0&\text{in $\mcw$},
}
where
\[\mcl_y:\phi\mapsto A\Delta_y\phi+B\cdot\nabla_y\phi\]
differentiates only in the variable $y\in\mcy$, and
\[\mcl_z:\phi\mapsto a\Delta_z\phi+b\cdot\nabla_z\phi\]
differentiates only in the variable $z\in\mcz$. Here and below, the boundary condition $\nu\cdot\nabla\varphi=0$ is required only when $\dr\mcw$ is nonempty. The coefficients $A,a\in\Er$, $B\in\Er^{n}$ and $b\in\Er^{p}$, for now, are assumed to be constant, while $c\in\mcc^{0,1}(\overline{\mcw})$ is Lipschitz. More general hypotheses on the coefficients will be stated below.

We are interested in the behaviour of the principal eigenvalue $k\in\Er$ and the principal eigenfunction $\varphi\in\mcc^{2}(\mcw)$ when the variable $y$ is very slow, in which case the elliptic operator ${\mcl_y+\mcl_z+c}$ becomes \enquote{strongly anisotropic}. Namely, for $\eps>0$, we define
\[\mcl^{\eps}_y\phi= \eps^2A\Delta_y\phi+\eps B\cdot\nabla_y\phi,\]
and we consider a solution $(\varphi_{\eps},k_{\eps})\in\mcc^2(\mcw)\times\Er$ of the principal eigenvalue problem:
\sys{\label{eq:intro_sf}
  \mcl^{\eps}_y\varphi_{\eps}+\mcl_z\varphi_{\eps}+c(y,z)\varphi_{\eps}&=k_{\eps}\varphi_{\eps}&\text{in $\mcw$},\\
  \nu\cdot\nabla\varphi_{\eps}&=0&\text{on $\dr\mcw$},\\
  \varphi_{\eps}&>0&\text{in $\mcw$}.
}
Our question is the following: How do $\varphi_{\eps}$ and $k_{\eps}$ behave as $\eps\to 0$? %A difficulty is raised because the operator $\mcl^{\eps}_y+\mcl_z+c$ is not assumed to be symmetric. To deal with this difficulty, we will use probabilistic methods.

For small $\eps>0$, the eigenvalue problem~\eqref{eq:intro_sf} is a two-scale system where the slow variable is $y$ and the fast variable is $z$. A way to get more insight into what this means is to state the following equivalent problem. With the change of variables $\varphi(y,z)=\widetilde{\varphi}\pth{\frac{y}{\eps},z}$, we obtain:
\sys{\label{eq:intro_sf2}
  \mcl_y\widetilde{\varphi}_{\eps}+\mcl_z\widetilde{\varphi}_{\eps}+c\pth{\eps y,z}\widetilde{\varphi}_{\eps}&=k_{\eps}\widetilde{\varphi}_{\eps}&\text{in $\mcw_{\eps}$},\\
  \nu\cdot\nabla\widetilde{\varphi}_{\eps}&=0&\text{on $\dr\mcw_{\eps}$},\\
  \widetilde{\varphi}_{\eps}&>0&\text{in $\mcw_{\eps}$},
}
where  $\mcw_{\eps}$ is defined by $\mcw_{\eps}:=\acc{\pth{\frac{y}{\eps},z}\in\Er^n\times\Er^p\tq (y,z)\in\mcw}$.
See Figure~\ref{fig:diff_scales} for the meaning of the two scales: in~\eqref{eq:intro_sf}, when we see the second-order terms as diffusive movements, the movements along $y$ are very small; in~\eqref{eq:intro_sf2}, the environment along $y$ is very wide.

\fig{7}{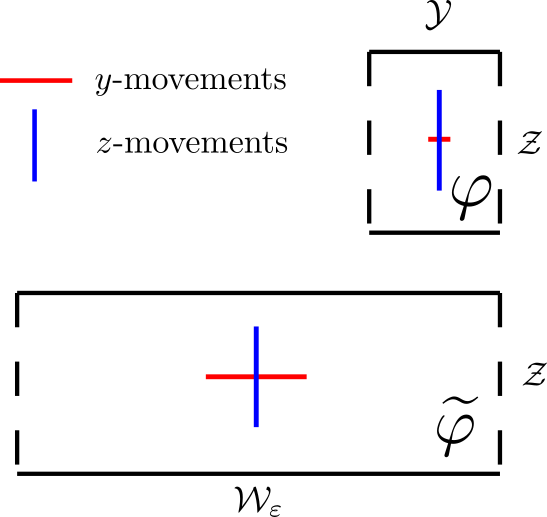}{\label{fig:diff_scales}The upper part corresponds to~\eqref{eq:intro_sf}; the lower part corresponds to~\eqref{eq:intro_sf2}. The crosses describe the ability of movements in either direction.} 

A first motivation for the study of the principal eigenvalue of an elliptic operator~$\mcl$ is its connection to the maximum principle for elliptic equations. If $\mcl$ acts on functions defined on a bounded and smooth domain $\mcw$, then the Krein-Rutman theory guarantees the existence and the uniqueness of the principal eigenvalue~$k(\mcl)$, and the maximum principle holds for the operator $\mcl$ if and only if $k(\mcl)<0$. Notions of principal eigenvalue, and their connections to the maximum principle, also exist when the domain is non-smooth~\cite{BNV94} or when the domain is unbounded~\cite{BerRos15}.
A second motivation arises in quantum mechanics. When $B\equiv 0$ and $b\equiv 0$ we get the Schrödinger operator; then $-c$ is a potential energy and the eigenvalues are energy levels for particles, and it is physically relevant to look for potentials optimising the eigenvalues, see \eg~\cite{Har84}.
Last, principal eigenvalues are often used to describe the long-time behaviour of a particular class of biological models. This is the main motivation for this work.
%The mathematical properties of the principal eigenvalue are thus interesting by themselves. Yet, our main motivation comes from the study of the long-time behaviour of biological models -- we will see indeed that their long-time behaviour is well described by principal eigenvalues.

\paragraph{Biological motivations.}
Consider the Fisher-KPP equation:
\[\dr_tu(t,y)=\Delta u+cu-u^2,\esp t\goq 0,\esp y\in\Er,\]
where $u(0,\cdot)\goq 0$ is nonzero and has a compact support.
Then $u(t,y)$ typically represents a density of individuals at time $t$ and position $y$, so that $u(t,\cdot)$ is the distribution of a population at time $t$. In this population dynamics setting, the constant $c>0$ is the intrinsic growth rate (\ie the growth rate of the population if the environment has infinite resources); and the $-u^2$ term stands for the competition between individuals. In their articles, Fisher~\cite{Fis37} and Kolmogorov, Petrovsky and Piskounov~\cite{KPP37} showed that the population eventually invades the environment at a constant velocity $v_{FKPP}=2\sqrt{c}$, in the following sense. If $v'>v_{FKPP}$, then $u(t,v't)\to 0$, and if $v''\in\ioo{0,v_{FKPP}}$, then $u(t,v''t)\to c>0$.

If the intrinsic growth rate $c$ depends on $y$ (say, if some places are more favourable than others), then we get the equation
\begin{equation}\label{eq:intro_evoleq_1}
  \dr_tu(t,y)=\Delta u+c(y)u-u^2,\esp t\goq 0,\esp y\in\Er.
\end{equation}
If the environment is too deleterious, \eg if $c(y)\ll 0$ for some $y\in\Er$, then the population may get extinct. Assume that $c$ is $1$-periodic. Then, a classical criterion to know whether the population persists or not depends on the sign of the principal eigenvalue of the symmetric operator
\[\mcl_{pers}:\phi\mapsto \Delta\phi+c(y)\phi\]
acting on $1$-periodic functions on $\Er$~\cite{BHR05}. Furthermore, Gärtner and Freidlin~\cite{GarFre79} showed the so-called Gärtner-Freidlin formula. Let $k(\lambda)$ be the principal eigenvalue associated with the nonsymmetric operator
\[\mcl_{GF}:\phi\mapsto \Delta\phi-2\lambda\dr_y\phi+(c(y)+\lambda^2)\phi\]
acting on $1$-periodic functions on $\Er$. If $k(0)>0$, then the population spreads at a velocity $v_{GF}$ given by
\[v_{GF}=\inf_{\lambda>0}\frac{k(\lambda)}{\lambda}.\]
This means that for all $v'>v_{GF}$, $u(t,v't)\to 0$ and for all $v''<v_{GF}$,
$u(t,x)$ converges locally uniformly, on $\abs{x}\loq v''t$, to
the unique $\phi>0$ satisfying $\Delta\phi+c\phi-\phi^2=0$.
%the profile $x\mapsto u(t,v''t+x)$ converges locally uniformly, up to a shift, to the unique $\phi>0$ satisfying $\Delta\phi+c\phi-\phi^2=0$.
The limiting profile $\phi$ is called the stationary state of the equation and, by uniqueness, is $1$-periodic.
In the study of the large-time behaviour of populations described by~\eqref{eq:intro_evoleq_1}, we thus need to study operators such as $\mcl_{pers}$ or $\mcl_{GF}$, which are similar to~\eqref{eq:intro_base}. 

If we want to study two distinct phenomena, or a two-dimensional space, we may add a new variable $z$ to the Fisher-KPP equation and get
\begin{equation}\label{eq:intro_evol_eq}
  \dr_tu(t,y,z)=\Delta_yu+\Delta_zu+c(y,z)u-u^2\esp t\goq 0,\esp y,z\in\Er.
\end{equation}
Then the persistence criterion and the Gärtner-Freidlin formula also hold in various cases~\cite{BHR05,BHR05-2,Wei02,BHN05}.

A typical setting is to consider $y$ as a phenotype and $z$ as a spatial position. Then $c(y,z)$ is the fitness of an individual of phenotype $y$ at position $z$, and~\eqref{eq:intro_evol_eq} becomes
\begin{equation}\label{eq:intro_evol_eq_nl}
  \dr_tu(t,y,z)=\Delta_yu+\Delta_zu+c(y,z)u-u\int_{\mcy}u(t,y',z)\de y'.
\end{equation}
The integral term comes from the fact that the competition occurs between all the individuals located at the same place, regardless of their phenotype. See~\cite{Pre04,ChaMel07} for the first references about this model. Again, recent works~\cite{ABR17,ACR13,BJS14} showed in various cases that the long-time behaviour of the solutions is governed by the principal eigenvalue of an elliptic operator.

Let us now explain the relevance of studying eigenvalue problems with a slow variable. In~\eqref{eq:intro_evol_eq_nl}, the diffusion in $y$ corresponds to mutations arising at each birth and death of an individual, and the diffusion in $z$ corresponds to spatial movements of the individuals. With this model in view, we understand that the diffusion in $y$ occurs on a much slower timescale than the diffusion in $z$. 
The initial motivation for this work was the study of the spreading speed of a population structured in space and phenotype with such a separation of timescales~\cite{b24}.
See also~\cite{BouMir15} for a study of this slow-fast dynamics from a Hamilton-Jacobi equation point of view.

Even when the variables $y$ and $z$ represent phenomena of the same nature, one may have to consider different scales. A first example occurs when $y$ and $z$ are spatial positions but on different scales of spatial heterogeneity: the variable $y$ represents, say, the heterogeneity within a large region, while the variable $z$ represents, say, the heterogeneity within a field.
A second example occurs when both variables $y$ and $z$ are the coordinates of a phenotype in a phenotypic landscape, but when mutations according to the $y$-direction are much rarer or much smaller than mutations according to the $z$-direction. See~\cite{HLMR20} for the study of such an anisotropic phenotypic landscape.

In all these cases, our results can be helpful to understand the interplay between the slow and the fast phenomenon, which are encoded in ~\eqref{eq:intro_sf} by the small parameter $\eps>0$.

We point out that different scaling limits (small diffusion and/or large transport term) have been studied in~\cite{BerHam02,CheLou08,CheLou12,HNR11}. Although they do not consider anisotropy, the setting in these works is similar, and some results in~\cite{HNR11,CheLou12} are reminiscent of ours. Finally, gradient estimates were given in~\cite{BerHam05} for the solutions of strongly anisotropic elliptic equations such as~\eqref{eq:intro_sf}.

\subsection{Main results}

\paragraph{General assumptions and notations.}
Let us describe our general assumptions. 
We define diffusion coefficients $A\in\mcc^{2,\alpha}(\overline{\mcy},\mcm_n(\Er))$ and $a\in\mcc^{2,\alpha}(\overline{\mcz},\mcm_p(\Er))$ (here $\mcm_q(\Er)$ is the set of $q\times q$-matrices), transport coefficients $B\in\mcc^{1,\alpha}(\overline{\mcy},\Er^{n})$ and $b\in\mcc^{1,\alpha}(\overline{\mcz},\Er^{p})$, and a Lipschitz-continuous potential $c\in\mcc^{0,1}(\overline{\mcw})$.
We assume that $A$ and $a$ are uniformly elliptic on $\overline{\mcw}$: there exist constants $\overline{A}>\underline{A}>0$ and $\overline{a}>\underline{a}>0$ such that for all $(y,z)\in\overline{\mcw}$, for all $\xi\in\Er^n$ and $\xi'\in\Er^p$,
\[\underline{A}\abs{\xi}^2\loq \xi\cdot (A(y,z)\xi)\loq\overline{A}\abs{\xi}^2,\esp \underline{a}\abs{\xi'}^2\loq \xi'\cdot (a(y,z)\xi')\loq\overline{a}\abs{\xi'}^2.\]
Last, we assume that $A$ and $a$ take their values in symmetric matrices; by the regularity of the coefficients, the solutions of the equations will be of class $\mcc^2$, so this symmetry assumption comes at no cost.

We let $\mcl^{\eps}_y$ be the operator defined by
\[\mcl^{\eps}_y:\phi\mapsto \eps^2\nabla_y\cdot(A\nabla_y\phi)+\eps\nabla_y\cdot\pth{\phi B},\]
and $\mcl_z$ be the operator defined by
\[\mcl_z:\phi\mapsto\nabla_z\cdot(a\nabla_z\phi)+\nabla_z\cdot\pth{\phi b}.\]
If $\dr\mcw\neq\emptyset$, we assume that the overall transport term is parallel to $\dr\mcw$:
\begin{equation}\label{eq:no_flux_dw_b}
  \pmat{B\\b}\cdot\nu=0\comment{on $\dr\mcw$}.
\end{equation}
For clarity, we write out here the variables of the coefficients:
\begin{equation}\label{eq:dependence_restriction}
  A(y)\vg a(z)\vg B(y)\vg b(z)\vg c(y,z).
\end{equation}
For $\eps>0$, we let $(\varphi_{\eps},k_{\eps})\in\mcc^{2,\alpha}(\mcw)\times\Er$ be the solution of the \emph{global} eigenvalue problem:
\sys{\label{eq:main}
  (\mcl^{\eps}_y+\mcl_z+c)\varphi_{\eps}&=k_{\eps}\varphi_\eps&\text{in $\mcw$},\\
  \varphi_{\eps}&>0&\text{in $\mcw$},\\
  \nu\cdot\mca\nabla\varphi_{\eps}&=0&\text{on $\dr\mcw$},\\
  \int_{\mcw}\varphi_{\eps}&=1,
}
where, for the statement of the boundary condition, we let $\mca:=\pmat{A&0\\0&a}\in\mcm_{n+p}(\Er)$ be the overall diffusion matrix.

\medskip

If $\dr\mcz$ is nonempty, we let $\nu_2$ be the outward unit normal on $\dr\mcz$.
For $y_0\in\overline{\mcy}$, we let $(\psi^{y_0},k^{y_0})\in\mcc^{2,\alpha}(\mcz)\times\Er$ be the solution of the \emph{local} eigenvalue problem:
\sys{\label{eq:local_eig_pb_general_thm}
  (\mcl_z+c(y_0,\cdot))\psi^{y_0}&=k^{y_0}\psi^{y_0}&\text{in $\mcz$},\\
  \psi^{y_0}&>0&\text{in $\mcz$},\\
  \nu_2\cdot a\nabla\psi^{y_0}&=0&\text{on $\dr\mcz$},\\
  \int_{\mcz}\psi^{y_0}&=1.
}
By the classical Krein-Rutman theory~\cite{KreRut48}, $(\varphi_{\eps},k_{\eps})$ and $(\psi^{y_0},k^{y_0})$ are well defined. 
If $\dr\mcy$ is nonempty, we extend the function $\varphi_{\eps}$ from $\mcy$ to $\overline{\mcy}$, \ie we consider $\varphi_{\eps}\in\mcc^{2,\alpha}(\overline{\mcy}\times\mcz)$. This is possible by virtue of the regularity assumptions on the coefficients and on $\dr\mcw$.

\paragraph{Behaviour of the principal eigenfunction.}

The first theorem means that the \enquote{slices} of the global principal eigenfunction converge, as $\eps\to 0$, to the local principal eigenfunctions. %As $\eps$ decreases, the coefficients $A,a,B,b,c$ are \enquote{more and more} independent of $y\in\mcy$, which suggests that the same holds true in the limit $\eps\to 0$.
We use the notation $\mcz'\subset\subset\mcz$ to say that $\overline{\mcz'}\subset\mcz$, and we set
\[\varphi_{\eps}(y_0,\mcz):=\int_{\mcz}\varphi_{\eps}(y_0,z)\de z.\]

\begin{thm}\label{thm:diffusion_zero} 
  Let $\beta\in\ioo{0,\alpha}$ and let $\mcz'\subset\subset\mcz$. For all $y_0\in\mcy$, as $\eps\to 0$, we have
  \begin{equation*}
    \frac{\varphi_\eps(y_0,\cdot)}{\varphi_\eps(y_0,\mcz)}\to \psi^{y_0}\esp\comment{in $\mcc^{2,\beta}(\mcz')$}. 
  \end{equation*}
  Moreover, the following convergence of measures holds uniformly in $y_0\in\overline{\mcy}$:
  \begin{equation*}
    \frac{\varphi_\eps(y_0,z)}{\varphi_\eps(y_0,\mcz)}\,\de z\to \psi^{y_0}(z)\de z\esp\comment{in total variation on $\mcz$},
  \end{equation*}
  namely:
  \begin{equation*}
    \sup_{y_0\in\overline{\mcy}}\ \cro{\sup_{A\subset\mcz,\,\text{$A$ measurable}}\ \abs{\int_{A}\pth{\frac{\varphi_\eps(y_0,z)}{\varphi_\eps(y_0,\mcz)}-\psi^{y_0}(z)}\de z}}\to 0.
  \end{equation*}
  
\end{thm}

Note that if $\mcz=\Te^p$, then we may take $\mcz'=\mcz$ and the convergence holds in $\mcc^{2,\beta}(\mcz)$.

 The result is straightforward when the coefficients $A$, $B$ and $c$ are independent of $y\in\mcy$: in this case, indeed, $\varphi_{\eps}=\psi^{y_0}$ is independent of $y_0$ and of $\eps$.
For the proof in the general case, we will use a probabilistic method. We will see $\varphi_{\eps}$ as the quasi-stationary distribution of a killed stochastic process $(Y^{\eps}_t,Z^{\eps}_t)_{t\goq 0}$. Namely, if $\tau$ is the death time of $(Y^{\eps}_t,Z^{\eps}_t)_{t\goq 0}$, we will have for all $E\subset\mcy$ and $F\subset\mcz$:
\[\lim_{t\to+\infty}\psf\pth{(Y^{\eps}_t,Z^{\eps}_t)\in E\times F\kt \tau>t}=\int_{E\times F}\varphi_{\eps}(y,z)\de y\de z.\]
As $\eps\to 0$, the process $Y^{\eps}_t$ will be slower and slower. This will allow us to forget the dependence in $Y^{\eps}_t$ and to get back to the simple case where the coefficients are independent of $y$. See Subsection~\ref{ss:definition_process} for the definition of the process $(Y^{\eps}_t,Z^{\eps}_t)_{t\goq 0}$.

We will see in the proof of Theorem~\ref{thm:diffusion_zero} why we need Assumption~\eqref{eq:dependence_restriction} on the coefficients.
  We will go by successive approximations from $(Y^{\eps}_t,Z^{\eps}_t)_t$, which corresponds to $\varphi_{\eps}$, to a stochastic process corresponding to $\psi^{y_0}$.
  For the approximations to be valid, we need a form of independence between $(Y^{\eps}_t)_t$ and $(Z^{\eps}_t)_t$;
  this independence is ensured by
  Assumption~\eqref{eq:dependence_restriction}.

Once Theorem~\ref{thm:diffusion_zero} is established, a natural question that arises is the behaviour of the principal eigenfunction in the slow direction.
  The family of principal eigenvalues $(k_{\eps})_{\eps}$ is bounded,
  so there is a subsequence that converges to some $k_0\in\Er$.
  The proofs of the following theorems suggest that
  if $\mcy$ has empty boundary, then for all $z_0\in\mcz$,
  there exists a measure~$\mu^{z_0}$ on the level set $\Gamma:=\acc{y\in\mcy\ /\ k^{y}=k_0}$,
  with total mass equal to $1$,
  such that for all $f\in\mcc^{\infty}(\mcy)$,
  along any subsequence such that $k_{\eps}\to k_0$,
  \[\int_{\mcy}f(y)\frac{\varphi_{\eps}(y,z_0)}{\varphi_{\eps}(\mcy,z_0)}\de y\to\int_{\Gamma} f(y)\mu^{z_0}(\de y).\]
  See Remark~\ref{req:conjecture_phi_slow} for more details.
  The proof of this conjecture and its extension to the case where $\mcy$ has a nonempty boundary are beyond the scope of this paper.

\paragraph{Behaviour of the principal eigenvalue.}
In two particular cases, we are able to find an exact relation between the limit of $k_{\eps}$ as $\eps\to 0$ and the local eigenvalues $k^{y}$, $y\in\overline{\mcy}$.
First, let us assume that $B\equiv 0$. 

\begin{thm}\label{thm:bzero}
  Assume $B\equiv 0$.
  Then $k_{\eps}$ converges as $\eps\to 0$ and
    \[\lim_{\eps\to 0}k_{\eps}=\max_{y\in\overline{\mcy}}k^{y}.\]
\end{thm}

If additionally to $B\equiv 0$, we assume that $b\equiv 0$, the operator $\mcl^{\eps}_y+\mcl_z+c$ becomes symmetric. The Rayleigh formula can then be applied directly, and the proof is simpler than in the general case. %The proof in this very particular case is given just before the proof in the general case (see Section~\ref{s:proof_cor}).

Theorem~\ref{thm:bzero} can be understood from a modelling perspective.
  When $B\equiv 0$ and $\eps\to0$, all \enquote{slices} of the environment, $c(y,\cdot)$, $y\in\overline{\mcy}$, are almost independent.
  Therefore, a species is able to survive in the environment if and only if there exists a favourable slice.
  This means that for small $\eps>0$, the condition $k_{\eps}>0$ amounts to the existence of $y_0\in\overline{\mcy}$ such that $k^{y_0}>0$.
  Thus we expect Theorem~\ref{thm:bzero} to hold.
  As we will see now, the independence between the slices fails when $B\not\equiv0$.

\medskip

Second, we relax the assumption that $B\equiv 0$, and assume instead that $\mcy$ is the {$1$-dimensional} torus. In particular, $B(y)\in\Er$ and $A(y)\in\Er$. We let
\[M:=\max_{y\in\mcy}\pth{k^y-\frac{B(y)^2}{4A(y)}}.\]
We let $j:\ifo{M,+\infty}\to\ifo{j(M),+\infty}$ be the bijection defined by
\[j:k\mapsto\int_0^1\sqrt{\frac{k-k^y+\frac{B(y)^2}{4A(y)}}{A(y)}}\de y.\]
Finally, we define
\[\gamma:=\int_0^1\frac{B(y)}{2A(y)}\de y.\]

\begin{thm}\label{thm:general} 
  Assume that $\mcy$ is the $1$-dimensional torus. %Let $k_{\eps}$ be the principal eigenvalue defined by~\eqref{eq:main}.
  Then $k_{\eps}$ converges as $\eps\to 0$ and:
  \begin{enumerate}
  \item If $\abs{\gamma}\loq j(M)$, then
    \[\lim_{\eps\to 0}k_{\eps}= M;\]
  \item If $\abs{\gamma}\goq j(M)$, then
    \[\lim_{\eps\to 0}k_{\eps}=j^{-1}\pth{\abs{\gamma}}.\]
  \end{enumerate}
\end{thm}

In the very specific case where all the coefficients are independent of the variable~$z$, Theorem~\ref{thm:general} is similar to Proposition 3.2 of~\cite{HNR11}.
If, moreover, $A=1$ and $B$ can be written in the form $B(y)=m'(y)$, then $\abs{\gamma}=0\loq j(M)$ and we get
\[\lim_{\eps\to 0}k_{\eps}=\max_{y\in\mcy}\pth{c(y)-\frac{(m'(y))^2}{4}}.\]
This is analogous to Theorem 1.3 of~\cite{CheLou12}, where $\mcy$ is assumed to be a smooth bounded domain in any dimension (contrarily to the setting of our Theorem~\ref{thm:general}).

The proof of Theorem~\ref{thm:general} is an adaptation of the proof of Proposition 3.2 of~\cite{HNR11}, but is more technical due to the interplay between the slow and fast variable.
  %Theorem~\ref{thm:diffusion_zero} allows us to forget the fast variable.
  As in~\cite{HNR11}, the proof hinges on the construction of an explicit viscosity solution to a Hamilton-Jacobi equation that arises naturally after a Hopf-Cole transform of the original eigenproblem.
  This explicit construction requires us to work in the $1$-dimensional setting;
  once the solution is constructed, we need to work in the periodic setting to deduce a result on the principal eigenvalue.

  Theorem~\ref{thm:general} shows that when there is a nonzero transport term, even in the particular case where~$\mcy$ is the $1$-dimensional torus, the 
  behaviour of the limiting principal eigenvalue is nontrivial.
  The study of the limit in the general case (bounded domain; higher dimension), therefore, is left to a future work.
Last, we point out that Theorem~\ref{thm:general} cannot be expected to hold in full generality if we drop Assumption~\eqref{eq:dependence_restriction}.
  Indeed, if~$A$ or~$B$ are allowed to depend on~$z$, then~$M$, $j$ and~$\gamma$ may not be well-defined.

\paragraph{Layout.} In Subsection~\ref{ss:definition_process}, we define a stochastic process which will be used in the proof of Theorem~\ref{thm:diffusion_zero} in Section~\ref{s:proof_thm}. Properties of this stochastic process are given in Section~\ref{s:qsd}. Finally, in Section~\ref{s:proof_cor}, we show Theorem~\ref{thm:bzero} and Theorem~\ref{thm:general}.

\subsection{Connection with a killed stochastic process}\label{ss:definition_process}

Let us define a stochastic process which will be at the core of the proof of Theorem~\ref{thm:diffusion_zero}.
We set
$c_m := \dabs{c}_{\infty}+\dabs{\nabla\cdot b}_{\infty}+\dabs{\nabla\cdot B}_{\infty}+1$
and
\[d(y,z):=c_m-c(y,z)>0.\]

\paragraph{Introduce a process $(\widetilde{Y}^{\eps}_t,\widetilde{Z}^{\eps}_t)_{t\goq 0}$.}

Let $m\in L^{1}(\mcw)$ satisfy
\[m\goq 0,\esp \int_{\mcw}m=1.\]
Let $\widetilde{u}^{\eps}$ be the solution of the following Cauchy problem:
\sys{\label{eq:pde_evol_yz}
  \dr_t\widetilde{u}^{\eps}(t,y,z)&=\eps^2\nabla_y\cdot(A\nabla_y\widetilde{u}^{\eps})+\eps \nabla_y\cdot(\widetilde{u}^{\eps}B)+\mcl_z\widetilde{u}^{\eps},&t>0\vg (y,z)\in\mcw,\\
  \nu\cdot\mca\nabla \widetilde{u}^{\eps}(t,y,z)&=0,&(y,z)\in\dr\mcw,\\
  \widetilde{u}^{\eps}(0,y,z)&=m(y,z),&(y,z)\in\mcw.}
In particular, thanks to the assumption~\eqref{eq:no_flux_dw_b}, we have for all $t\goq 0$,
\[\int_{\mcw}\widetilde{u}^{\eps}(t,\cdot)=1.\]
Using Theorem 7.1 of~\cite{SatUen65}, we construct a Markov process $(\widetilde{Y}^{\eps}_t,\widetilde{Z}^{\eps}_t)_{t\goq 0}$, with $\widetilde{Y}^{\eps}_t\in\mcy$ and $\widetilde{Z}^{\eps}_t\in\mcz$, with continuous paths, and which has a density at time $t$ given by $\widetilde{u}^{\eps}(t,\cdot,\cdot)$. %Moreover $(\widetilde{Y}^{\eps}_t,\widetilde{Z}^{\eps}_t)_{t\goq 0}$ has continuous paths almost surely. { add details... See Itô [18] for why their operator is strange}
Typically, if $a=\sigma^2$ and $A=\Sigma^2$ are constant, we should think of the process $(\widetilde{Y}^{\eps}_t,\widetilde{Z}^{\eps}_t)_{t\goq 0}$ as a solution of the SDE:
\sys{\label{eq:sde}
  &\de\widetilde{Y}^{\eps}_t=-\eps B(\widetilde{Y}^{\eps})\de t+\sqrt{2}\eps \,\Sigma\de \xi^Y_t,\\
  &\de\widetilde{Z}^{\eps}_t=-b(\widetilde{Z}^{\eps})\de t+\sqrt{2}\,\sigma\de \xi^Z_t,\\
  &\text{reflection on $\dr\mcw$},
}
with initial distribution $(\widetilde{Y}^{\eps}_0,\widetilde{Z}^{\eps}_0)\sim m$,
  and where $\xi^Y$ and $\xi^Z$ are independent $n$-dimensional and $p$-dimensional Brownian motions respectively.
%Typically, we should think of the process $(\widetilde{Y}^{\eps},\widetilde{Z}^{\eps})$ as satisfying the SDEs:
%\sys{\label{eq:sde}
%  &\de\widetilde{Y}^{\eps}_t=-\pth{\eps B(\widetilde{Y}^{\eps})-\eps^2\nabla_yA(\widetilde{Y}^{\eps})}\de t+\sqrt{2}\eps A(\widetilde{Y}^{\eps})\de W^Y_t,\\
%  &\de\widetilde{Z}^{\eps}_t=-\pth{b(\widetilde{Z}^{\eps})-\nabla_za(\widetilde{Z}^{\eps})}\de t+\sqrt{2}a(\widetilde{Z}^{\eps})\de W^Z_t,\\
%  &\text{reflection on $\dr\mcw$},
%  }
%  where $W^Y$ and $W^Z$ are independent $n$-dimensional and $p$-dimensional Brownian motions respectively.

\paragraph{Kill the process $(\widetilde{Y}^{\eps}_t,\widetilde{Z}^{\eps}_t)_{t\goq 0}$.}

Now, we construct a killed process $({Y}^{\eps}_t,{Z}^{\eps}_t)$ based on $(\widetilde{Y}^{\eps}_t,\widetilde{Z}^{\eps}_t)$. We first need a definition, which will be used several times throughout this work. The definition is illustrated by Figure~\ref{fig:poisson}.
\begin{dfi}\label{dfi:poisson}
Let $(D_t)_{t\goq 0}$ be a stochastic process with nonnegative values. Let $\Pi$ be a Poisson point process on $\Er_+\times\Er_+$ with intensity $1$, independent of $(D_t)_{t\goq 0}$. Define
\[\tau((D_t)_{t\goq 0}):=\inf\acc{t\goq 0\tq \Pi\text{ has a point in $\acc{(s,y)\tq 0\loq s\loq t\vg 0\loq y\loq D_t}$}}.\]
We say that $\tau((D_t)_{t\goq 0})$ is a clock with rate $(D_t)_{t\goq 0}$.
\end{dfi}

\fig{7}{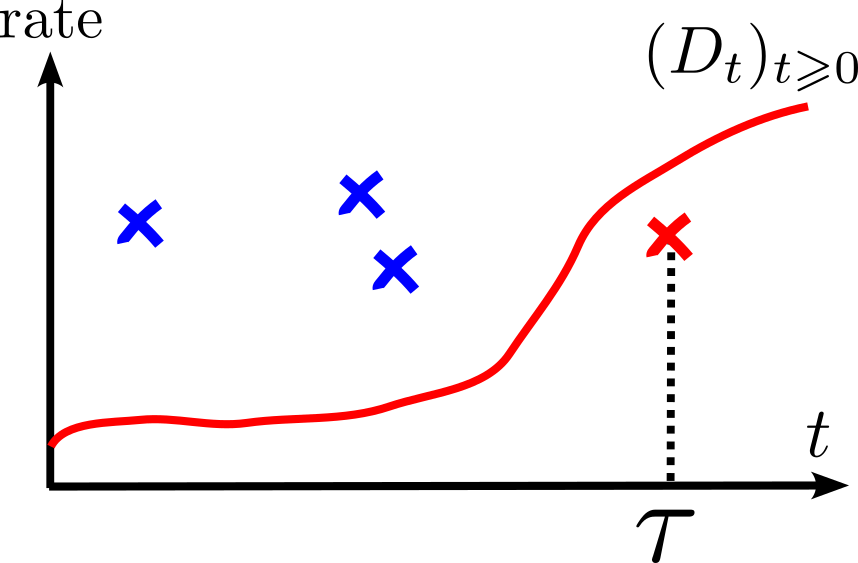}{\label{fig:poisson}The crosses are the points of $\Pi$. We are interested in $\tau((D_t)_{t\goq 0})$, the first time a cross is under the curve.}

Using this definition, we let $\tau^{\eps}$ be a clock with rate $(d(\widetilde{Y}^{\eps}_t,\widetilde{Z}^{\eps}_t))_{t\goq 0}$.
We let $\partial\notin\mcw$ be a \enquote{cemetery state}. We now define the killed process $(Y^{\eps}_t,Z^{\eps}_t)_{t\goq 0}$ by:
\[(Y^{\eps}_t,Z^{\eps}_t)=\left\{
\begin{aligned}
  &(\widetilde{Y}^{\eps}_t,\widetilde{Z}^{\eps}_t),&t<\tau,\\
  &\dr,&t\goq\tau.
\end{aligned}
\right.\]
We will see ${\tau}^{\eps}$ as a \enquote{death time} for the particle, and work conditionally on $\acc{\tau^{\eps}>t}$. See Figure~\ref{fig:desc_ytzt}, left, for a typical trajectory before a death occurs; and see Figure~\ref{fig:desc_ytzt}, right, for the effect of the conditioning on $\acc{\tau^{\eps}>t}$.

\begin{figure}[h] 
  \begin{center}
    \includegraphics[height=4cm]{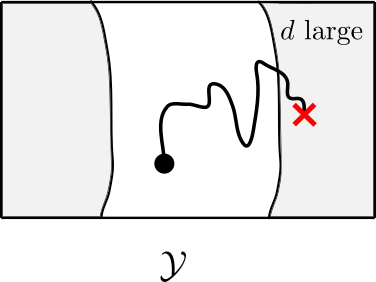}\esp    
    \includegraphics[height=4cm]{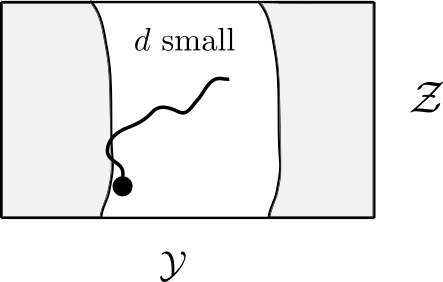}
  \end{center}

  \caption{\label{fig:desc_ytzt}\emph{Grey zone:} high death rate. \emph{White zone:} low death rate. \emph{Left:} The process $(Y^{\eps}_t,Z^{\eps}_t)_{t\goq 0}$ dies at a high rate when it is in the grey zone; \emph{Right:} Conditioned on $\acc{\tau^{\eps}>t}$, the process $(Y^{\eps}_t,Z^{\eps}_t)_{t\goq 0}$ avoids the grey zone.}
  
\end{figure}

We denote by $\psf_m$ the joint law of $(\tau,\, (\widetilde{Y}^{\eps}_t,\widetilde{Z}^{\eps}_t)_{t\goq0})$ when the initial condition in~\eqref{eq:pde_evol_yz} is~$m$. The following theorem implies that $\varphi_{\eps}$ is the unique quasi-stationary distribution for the killed process $(Y^{\eps}_t,Z^{\eps}_t)_{t\goq 0}$. We will use it to find properties of $\varphi_{\eps}$ by studying the process $(Y^{\eps}_t,Z^{\eps}_t)_{t\goq 0}$.
Let $\dabs{\cdot}_{TV}$ denote the total variation norm for (signed) measures on~$\mcw=\mcy\times\mcz$. In particular, for $f:\mcw\to\Er$,
\[\dabs{f(x)\de x}_{TV}:=\sup_{A\subset\mcw,\, \text{$A$ measurable}}\ \abs{\int_{\mcz}f}.\]
%  \[\dabs{f(x)\de x}_{TV}:=\sup_{A\subset\mcw,\, \text{$A$ measurable}}\ \int_{A}|f|.\]

\begin{thm}\label{thm:cv_qsd_yz}
    For all $\eps>0$, there exist $C,\chi>0$ such that for all initial distribution $m$ on~$\mcw$, and for all $t>0$,
    \begin{align*}
      \dabs{\varphi_{\eps}(y,z)\de y\de z-\psf_{m}\pth{(Y^{\eps}_t,Z^{\eps}_t)\in \de y\de z\kt \tau>t}}_{TV}<Ce^{-\chi t}.
%    \dabs{\varphi_{\eps}(\cdot)-\psf_{m}\pth{(Y^{\eps}_t,Z^{\eps}_t)\in \cdot\kt \tau>t}}_{TV}<Ce^{-\chi t}.
  \end{align*}
\end{thm}

In the next section, we will prove a more general version of Theorem~\ref{thm:cv_qsd_yz}.

\section{Results about quasi-stationary distributions}\label{s:qsd}

\subsection{General setting}

As we will need an equivalent of Theorem~\ref{thm:cv_qsd_yz} in another setting, we prove a more general version. We proceed to a construction similar to that of $(Y^{\eps},Z^{\eps})$, but with simpler notations.

Let $\mcy$ and $\mcz$ be defined as in the introduction and let, for this section, $\mcw:=\mcy\times\mcz$ or $\mcw:=\mcz$.
Let~$q$ be the dimension of~$\mcw$ (thus either $q=n+p$ or $q=p$).
%=\mcy\times\mcz$ where $\mcy$ is a torus or a smooth domain of $\Er^{n_1}$, and  where $\mcy$ is a torus or a smooth domain of $\Er^{n_2}$.
Again, we let $\nu$ be the outward unit normal on~$\dr\mcw$.
Throughout Section~\ref{s:qsd}, we do not treat the variables~$y$ and~$z$ as separated variables, and we let $w\in\mcw$ be the generic element of~$\mcw$.
Let $\mca\in\mcc^{2,\alpha}(\overline{\mcw},\mcm_q(\Er))$ and $\mcb\in\mcc^{1,\alpha}(\overline{\mcw},\Er^{q})$.
Let $\Lambda$ be the operator defined by
\[\Lambda u=\nabla\cdot(\mca \nabla u)+\nabla\cdot(u\mcb).\]
When $\mcw=\mcy\times\mcz$, we may think of $\Lambda$ as the operator unifying the $\mcy$- and the $\mcz$-variables, that is:
\[\mca=\pmat{A&0\\0&a},\esp\mcb=\pmat{B\\b}.\]
We assume that the transport term is parallel to $\dr\mcw$:
\begin{equation}\label{eq:no_flux_dw_beta}
  \mcb\cdot\nu=0\comment{on $\dr\mcw$}.
\end{equation}
This assumption is a translation of~\eqref{eq:no_flux_dw_b} in this new framework.
We further assume that $\mca$ is uniformly elliptic on $\overline{\mcw}$: there exist $\overline{\mca}>\underline{\mca}>0$ such that for all $w\in\overline{\mcw}$ and $\xi\in\Er^q$,
\[\underline{\mca}\abs{\xi}^2\loq \xi\cdot \mca(w)\xi\loq\overline{\mca}\abs{\xi}^2.\]
Last, we assume that $\mca$ takes its values in symmetric matrices. %; by the regularity of the coefficients, the solutions of the equations will be of class $\mcc^2$, so this symmetry assumption comes at no cost.}

Let $m\in L^{1}(\mcw)$ satisfy
\[m\goq 0,\esp\int_{\mcw}m=1.\]
Let $\widetilde{u}$ be the solution of the following Cauchy problem:
\sys{\label{eq:cp_evol_wtilde}
  \dr_t\widetilde{u}(t,w)&=\Lambda \widetilde{u},&t>0\vg w\in\mcw,\\
  \nu\cdot\mca\nabla \widetilde{u}(t,w)&=0,&w\in\dr\mcw,\\
  \widetilde{u}(0,w)&=m(w),&w\in\mcw.}
In particular, thanks to the assumption~\eqref{eq:no_flux_dw_beta}, we have for all $t\goq 0$,
\[\int_{\mcw}\widetilde{u}(t,\cdot)=1.\]
Using Theorem 7.1 of~\cite{SatUen65}, we construct a Markov process $(\widetilde{W}_t)_{t\goq 0}$, with continuous paths, and which has a density at time $t$ given by $\widetilde{u}$.
If $\mca=\sigma^2$ is constant, we should think of the process $(\widetilde{W}_t)_{t\goq 0}$ as the solution to the SDE
  \[\de \widetilde{W}_t=-\mcb(\widetilde{W}_t)\de t+\sqrt{2}\sigma\de \xi_t,\]
  with reflection on~$\dr\mcw$ and with initial distribution~$m$,
  and
  where $\xi_t$ is a standard Brownian motion in~$\Er^q$.
Let
\[\mce:=\acc{f\in\mcc^2(\overline{\mcw})\tq \nu\cdot\mca \nabla f=0\text{ on $\dr\mcw$}}.\]
  Let $\Lambda^*$ be the adjoint of $\Lambda$, defined on $\mce$ by $\Lambda^*:f\in\mce\mapsto \nabla\cdot(\mca\nabla f)-\mcb\cdot\nabla f$.
  The process $(\widetilde{W}_t)_{t\goq 0}$ has a semigroup on $\mcc^0(\mcw)$
  with infinitesimal generator $\Lambda^*$ defined on the domain $\mce$.
  This means that for all $f\in\mce$,
%    \[\esf_w\cro{f(\widetilde{W}_t)}\to f(w) \qquad\text{as $t\to0$},\]
%    and
    \[\frac{\esf_{\delta_w}\cro{f(\widetilde{W}_t)}-f(w)}{t}\to \Lambda^*f(w) \qquad\text{as $t\to0$},\]
    uniformly in $w\in\mcw$, see Chapter~1 of~\cite{EthKur86}. Here, $\esf_{\delta_w}$ corresponds to the expectation when $\widetilde{W}_0=w$.

Now, we construct a killed process $(W_t)_{t\goq 0}$ based on $(\widetilde{W}_t)_{t\goq 0}$. We let $d\in\mcc^{0,1}(\overline{\mcw})$ be a Lipschitz-continuous function such that $\displaystyle\inf_{\mcw}d\goq 1+\dabs{\nabla\cdot\mcb}_{\infty}$. We let $\tau$ be a clock with rate $(d(\widetilde{W}_t))_{t\goq 0}$ (in the sense of Definition~\ref{dfi:poisson}).
We let $\partial\notin\mcw$ be a \enquote{cemetery state}.
We define the killed process $(W_t)_{t\goq 0}$ by:
\[W_t:=
\left\{
\begin{aligned}
  &\widetilde{W}_t&t<\tau,\\
  &\dr&t\goq \tau.
\end{aligned}
\right.
\]
We denote by $\psf_m$ the joint law of $(\tau, (\widetilde{W}_t)_{t\goq 0})$ when the initial distribution of $(\widetilde{W}_t)_t$ is $m$, and by $\esf_m$ the corresponding expectation. 

Now, we let $(\Phi,k)\in\mcc^{2,\alpha}(\mcw)\times\Er$ satisfy the following principal eigenvalue problem:
\syss{(\Lambda-d(w))\Phi(w)&=k\Phi(w),&w\in\mcw,\\
  \Phi(w)&>0,&w\in\mcw,\\
\nu\cdot\mca\nabla\Phi(w)&=0,&w\in\dr\mcw.
}
We normalise $\Phi$ so that \[\int_{\mcw}\Phi=1.\]
The following theorem shows that the principal eigenfunction $\Phi$ is a quasi-stationary distribution for the killed process $(W_t)_{t\goq 0}$. 

\begin{thm}\label{thm:cv_qsd}
  There exist $C,\chi>0$ such that for all initial distribution $m$ on $\mcw$, and for all $t>0$,
  \begin{align*}
    \dabs{\Phi(w)\de w-\psf_{m}\pth{W_t\in \de w\kt \tau>t}}_{TV}<Ce^{-\chi t}.
%    \dabs{\Phi(\cdot)-\psf_{m}\pth{W_t\in \cdot\kt \tau>t}}_{TV}<Ce^{-\chi t}.
  \end{align*}
  Let $\mcb_M\goq \dabs{\mcb}_{\infty}$. We may choose the constants $C$ and $\chi$ so that they only depend on $\underline{\mca}$, $\overline{\mca}$, $\mcb_M$, $\dabs{d}_{\infty}$ and $\mcw$.
\end{thm}

Theorem~\ref{thm:cv_qsd_yz} is a consequence of Theorem~\ref{thm:cv_qsd}. In Subsections~\ref{ss:pty_killed} and~\ref{ss:proof_thm_qsd}, we will prove Theorem~\ref{thm:cv_qsd}.

\subsection{Properties of the killed stochastic process}\label{ss:pty_killed}

Let us now study the law of $(W_t)_{t\goq 0}$.
We first give the generator of the process; after that, we will give its density.

We denote by $\mce^{\dr}$ the set of functions $f:\mcw\cup\acc{\dr}\to\Er$ such that $f_{|\mcw}\in\mce$.

\begin{lem}\label{lem:generator_xtst}

  The Markov process $(W_t)_{t\goq 0}$ has a semigroup on $\mcc^0(\mcw)$
  with infinitesimal generator ${\mcg}$ defined on the domain $\mce^{\dr}$ by:
  \begin{align*}
    \mathcal{G}f(w)&=(\Lambda^*-d(w)) f(w)+d(w)f(\partial),\esp w\in\mcw,\\
    \mathcal{G}f(\partial)&=0.
  \end{align*}
  In particular, %Moreover,
  for all $f\in\mce^{\dr}$, the convergence of
  \[  
  \frac{1}{h}\pth{\esf_{\delta_w}\cro{f(W_h)}-f(w)}
  \]
  to $\mcg f(w)$ as $h\to 0$ is uniform in $w\in\mcw$.
\end{lem}

\begin{proof}
  Since $\partial$ is a cemetery state, we have $\mathcal{G}f(\partial)=0$. Now, let us compute $\mcg f(w)$ for $w\in\mcw$.
  For $h>0$,
  \begin{align}
    \esf_{\delta_w}\cro{f(W_h)}&=\esf_{\delta_w}\cro{f(W_h)\ind_{{\tau>h}}}+\esf_{\delta_w}\cro{f(\partial)\ind_{{\tau\loq h}}}\nonumber\\
    &=\esf_{\delta_w}\cro{f(W_h)\ind_{{\tau>h}}}+f(\partial)\psf_{\delta_w}\pth{\tau\loq h}.\label{eq:decomposition_fwh}
  \end{align}
  First, since $d$ is continuous on $\overline{\mcw}$, $d$ is also uniformly continuous on $\overline{\mcw}$; thus
  \[\psf_{\delta_w}(\tau\loq h)=hd(w)+o_{h\to 0}(h),\]
  and the $o_{h\to 0}(h)$ is uniform in $w$.
  Second,
  \begin{align*}
    \esf_{\delta_w}\cro{f(W_h)\ind_{\tau> h}}&=\esf_{\delta_w}\cro{f(\widetilde{W}_h)\ind_{\tau> h}}\\
    &=\esf_{\delta_w}\cro{f(\widetilde{W}_h)}-\esf_{\delta_w}\cro{(f(\widetilde{W}_h)
      -f(w))\ind_{\tau\loq h}}-\esf_{\delta_w}\cro{f(w)\ind_{\tau\loq h}}.
  \end{align*}
  We have $\esf_{\delta_w}\cro{(f(\widetilde{W}_h)-f(w))\ind_{\tau\loq h}}=O_{h\to0}(h^2)$
    and $\esf_{\delta_w}\cro{f(w)\ind_{\tau\loq h}}=f(w)\psf_{\delta_w}(\tau\loq h)$.
    Therefore, using the fact that $\Lambda^*$ is the infinitesimal generator of the process $(\widetilde{W}_t)_t$, we have
  \begin{align*}
    \esf_{\delta_w}\cro{f(W_h)\ind_{\tau> h}}&=f(w)+h\Lambda^*f(w)-hd(w)f(w)+o_{h\to 0}(h).
  \end{align*}
  We get
  \begin{align*}
    \esf_{\delta_w}\cro{f(W_h)}&=
    f(w)+h(\Lambda^*-d)f(w)+hd(w)f(\partial)+o_{h\to 0}(h).
  \end{align*}
  Moreover, the $o_{h\to 0}(h)$ term remains uniform in $w\in\mcw$. Therefore,
  \[\mcg f=(\Lambda^*-d(w))f(w)-d(w)f(\partial),\]
  and the convergence to $\mcg f$ is uniform on $\mcw$. 
\end{proof}

We are now ready to give the density of $(W_t)_{t\goq 0}$. 
Consider the evolution problem
\sys{\label{eq:cp_density_w}
  \dr_tu(t,w)&=(\Lambda-d(w))u(t,w),&t>0\vg w\in\mcw,\\
  u(0,\cdot)&=m,&w\in\mcw,\\
  \nu\cdot\mca\nabla u(w)&=0,&w\in\dr\mcw.
}

\begin{lem}\label{lem:density_pde} 
The solution $u(t,\cdot)$ of~\eqref{eq:cp_density_w} is the density of the law of $W_t$. Namely, for all $t\goq 0$, for all measurable bounded function ${f:\mcw\cup\acc{\partial}\to\Er}$:
  \begin{equation}\label{eq:u_density_w}
    \esf_{m}\cro{f(W_t)}=\psf_{m}\pth{\tau<t}f(\partial)+\int_{\mcw}f(w)u(t,w)\de w.
  \end{equation}
\end{lem}

\begin{proof}
  The law of $W_t$ restricted to $\mcw$ is absolutely continuous with respect to the law of~$\widetilde{W}_t$, which has a density with respect to the Lebesgue measure. For $t>0$, we let $\zeta(t,\cdot)$ be the density of the law of $W_t$ (restricted to $\mcw$) with respect to the Lebesgue measure.
  We show that $\zeta$ satisfies~\eqref{eq:cp_density_w} in a weak sense.

Let $u$ be the solution of the Cauchy Problem~\eqref{eq:cp_density_w}. Let $T>0$ and $F\in\mce$. We have
\begin{align*}
  \int_{\mcw}(u(T,w)-u(0,w))F(w)\de w&=\int_0^T\int_{\mcw}F(w)\dr_tu(t,w)\de w\de t\\
  &=\int_0^T\int_{\mcw}F(w)\cro{(\Lambda-d)u(t,w)}\de w\de t.
\end{align*}
Using the divergence theorem, we obtain:
  \begin{equation}\label{eq:weak_form_density}
    \forall T>0\vg\forall F\in\mce,\esp\int_{\mcw}(u(T,w)-u(0,w)) F(w)\de w=\int_0^T\int_{\mcw}u(t,w) [(\Lambda^*-d)F(w)]\de w\de t.
  \end{equation}
  We call~\eqref{eq:weak_form_density} the weak form of the Cauchy problem~\eqref{eq:cp_density_w}.

  We now show that the density $\zeta$ satisfies~\eqref{eq:weak_form_density}. 
   Let $F\in\mce$ and set $F(\dr)=0$, so that we can also see $F$ as an element of $\mce^{\dr}$. Let $\mathscr{F}_t=\sigma((W_s),0\loq s\loq t)$ be the $\sigma$-field generated by $(W_s)_{s\goq 0}$. Then
  \[F(W_t)-F(W_0)-\int_0^t(\Lambda^*-d)F(W_s)\de s\]
  is a $\mathscr{F}_t$-martingale, see \eg \cite[Proposition 4.1.7]{EthKur86}. 
  Thus for all $T>0$,
  \begin{align*}
    \Esp{F(W_T)-F(W_0)}&=\int_0^T\int_{\mcw}\zeta(t,w)[(\Lambda^*-d) F(w)]\de w\de t.
  \end{align*}
  We obtain:
  \begin{equation*}
    \forall T>0\vg\forall F\in\mce,\esp\int_{\mcw}\pth{\zeta(T,w)-\zeta(0,w)}F(w)\de w=\int_0^T\int_{\mcw}\zeta(t,w) [(\Lambda^*-d)F(w)]\de w\de t.
  \end{equation*}
  This means precisely that~$\zeta$ satisfies~\eqref{eq:weak_form_density}.
  It is standard that the weak formulation~\eqref{eq:weak_form_density} has a unique solution with initial condition $m$. Therefore, $u=\xi$.
\end{proof}

\subsection{Proof of Theorem~\ref{thm:cv_qsd}}\label{ss:proof_thm_qsd}

Before proving Theorem~\ref{thm:cv_qsd}, we give a preliminary result which implies that the principal eigenfunction~$\Phi$ is a quasi-stationary distribution for the process $(W_t)_{t\goq 0}$.

\begin{ppn}\label{ppn:qsd}
  For all measurable set $B\subset\mcw$, and for all $t>0$,
\begin{align*}
  \Phi(B)&=\psf_{\Phi}\pth{W_t\in B\kt \tau>t}.
\end{align*}  
\end{ppn}

\begin{proof}
  Starting from $m=\Phi$, the solution of the system~\eqref{eq:cp_density_w} is:
  \[u(t,w)=e^{kt}\Phi(w).\]
  We point out that, due to the positivity of $d$, we have $k\loq 0$; thus $e^{kt}\in(0,1]$.
  We then apply Lemma~\ref{lem:density_pde} with initial condition~$\Phi$ and with the test function $f := \ind_{B}$ with $f(\dr)=0$. We find
  \[\psf_{\Phi}\pth{W_t\in B}=\int_Be^{kt}\Phi(w)\de w=e^{kt}\psf_{\Phi}\pth{W_0\in B}.\]
  Now, we note that $\psf_{\Phi}\pth{\tau>t}=\psf_{\Phi}\pth{W_t\in\mcw}$.
  Applying again Lemma~\ref{lem:density_pde} with the test function $f=\ind_{\mcw}$, we conclude that
  \[\psf_{\Phi}\pth{\tau>t}=e^{k t}.\]
  Therefore,
  \begin{align*}
    \psf_{\Phi}\pth{W_t\in B\kt \tau>t}&=\frac{\psf_{\Phi}\pth{W_t\in B, \tau>t}}{\prb{\tau>t}}\\
    &=\frac{\psf_{\Phi}\pth{W_t\in B}}{\prb{\tau>t}}\\
    &=\psf_{\Phi}\pth{W_0\in B}.
  \end{align*}
  This is the statement of the proposition. 
\end{proof}

We are now ready to prove Theorem~\ref{thm:cv_qsd}. % that $\pmat{\varphi\\\psi}$ is a quasi-stationary distribution.

\begin{proof}[Proof of Theorem~\ref{thm:cv_qsd}]
  We apply Theorem 1.1 of~\cite{ChaVil14}. 
  It is enough to check that there exists a probability measure $\mu$ on $\mcw$ such that the following two assumptions hold:
  \begin{enumerate}
  \item \textbf{Assumption A1.} There exist $t_0>0$ and $c_1>0$ such that for all $w\in\mcw$, for all measurable $U\subset\mcw$,
    \begin{equation}\label{eq:ass_a1_chavil}
      \psf_{\delta_w}\pth{W_{t_0}\in U\kt \tau>t_0}\goq c_1\mu(U);
    \end{equation}
  \item \textbf{Assumption A2.} There exists $c_2>0$ such that for all $t>0$, for all $w\in\mcw$,
  \begin{equation}\label{eq:ass_a2_chavil}
    \psf_{\delta_w}\pth{\tau>t}\loq c_2\psf_{\mu}\pth{\tau>t}.
  \end{equation}
  \end{enumerate}
  Here $\psf_{\delta_w}$ denotes the law corresponding to the the fundamental solution $\widetilde{u}^{(w)}(t,w')$ associated with~\eqref{eq:cp_evol_wtilde}, that is:
  \syss{
  \dr_t\widetilde{u}^{(w)}(t,w')&=\Lambda \widetilde{u}^{(w)},&t>0\vg w'\in\mcw,\\
  \nu\cdot\mca\nabla \widetilde{u}^{(w)}(t,w')&=0,&w'\in\dr\mcw,\\
  \lim_{t\to 0}\widetilde{u}^{(w)}(t,\cdot)&=\delta_{w}.}
    By Proposition~\ref{ppn:qsd}, $\Phi$ is a quasi-stationary distribution of $(W_t)_{t\goq 0}$. If Assumptions A1 and A2 hold, then $\Phi$ must be the unique quasi-stationary distribution of $(W_t)_{t\goq 0}$ and the exponential convergence to $\Phi$ stated in Theorem~\ref{thm:cv_qsd} must hold. 

    \paragraph{Step 1. Check that Assumption A1 holds.}
  We let $\mu:=\frac{\de w}{\abs{\mcw}}$ be the uniform law on $\mcw$.
  We first check that~\eqref{eq:ass_a1_chavil} holds with $t_0:=2$ and for $c_1>0$ small enough. 
  Take $\mcb_M\goq\dabs{\mcb}_{\infty}$. Using the regularity of $\mcw$, we extend the solution $\widetilde{u}^{(w)}$ by orthogonal reflection with respect to $\dr\mcw$. We get a function defined on a set $\widetilde{\mcw}$ containing $\overline{\mcw}$. The operator $\Lambda$ is extended to a uniformly elliptic operator $\widetilde{\Lambda}$ on $\widetilde{\mcw}$. This allows one to apply the parabolic Harnack inequality on the whole $\mcw$: there exists $C>0
  $ depending only on $\underline{\mca}$, $\overline{\mca}$, $\mcb_M$ and~$\mcw$ such that for all $t> 0$,
  \[\inf_{w'\in\mcw}\widetilde{u}^{(w)}(t+1,w')\goq C\sup_{w'\in\mcw}\widetilde{u}^{(w)}(t,w').\]
  Since we have also, for all $t>0$,
  \[\int_{\mcw}\widetilde{u}^{(w)}(t,w')\de w'=1,\]
  we obtain: for all $t> 0$,
  \[\inf_{w'\in\mcw}\widetilde{u}^{(w)}(t+1,w')\goq \eta_1:=\frac{C}{\abs{\mcw}}.\]
%  By the strong maximum principle (and the Neumann boundary condition), we have \[\inf_{w'\in\mcw}\widetilde{u}^{(w)}(1,w')>0.\] 
%  Since $(w,w')\mapsto \widetilde{u}^{(w)}(1,w')$ is continuous, we obtain:
%  \[\eta_1:=\inf_{w,w'\in\mcw}\widetilde{u}^{(w)}(1,w')>0.\]
  We let $\underline{\tau}$ be a clock independent of $(W_t)_{t\goq 0}$, with rate $D:=\dabs{d}_{\infty}$, and coupled with $\tau$ so that $\underline{\tau}\loq \tau$ almost surely (this is possible by defining $\tau$ and $\underline{\tau}$ with the same Poisson point process taken independent of $(W_t)_{t\goq 0}$). We obtain:
  \begin{align*}
    \psf_{\delta_w}\pth{W_2\in U\kt \tau>2}&\goq  \psf_{\delta_w}\pth{W_2\in U\vg \underline{\tau}>2}\\
    &= \psf_{\delta_w}\pth{W_2\in U}e^{-2D}\\
    &=e^{-2D}\int_U\widetilde{u}^{(w)}(2,w')\de w'\\
    &\goq\eta_1 e^{-2D}\abs{U}=\eta_1 e^{-2D}\abs{\mcw}\mu(U).
  \end{align*}
  Assumption A1 is satisfied with $c_1:=\eta_1 e^{-2D}\abs{\mcw}$.

  \paragraph{Step 2. Check that Assumption A2 holds.}
   We let ${u}^{(w)}(t,w')$ be the fundamental solution associated with~\eqref{eq:cp_density_w}, that is:
  \syss{
  \dr_t{u}^{(w)}(t,w')&=(\Lambda-d(w)) {u}^{(w)},&t>0\vg w'\in\mcw,\\
  \nu\cdot\mca\nabla {u}^{(w)}(t,w')&=0,&w'\in\dr\mcw,\\
                 \lim_{t\to 0}{u}^{(w)}(t,\cdot)&=\delta_{w}.}
    Again, we can apply the parabolic Harnack inequality on the whole $\mcw$: there exists $C'>0$ depending only on $\underline{\mca}$, $\overline{\mca}$, $\mcb_M$, $D:=\dabs{d}_{\infty}$ and~$\mcw$ such that for all $t> 0$,
    \[\inf_{w'\in\mcw}{u}^{(w)}(t+1,w')\goq C'\sup_{w'\in\mcw}{u}^{(w)}(t,w').\]
  Since we have also, for all $t>0$,
  \[\int_{\mcw}{u}^{(w)}(t,w')\de w'\goq e^{-Dt},\]
  we obtain in particular: 
  \[\inf_{w'\in\mcw}{u}^{(w)}(2,w')\goq \eta_2\]
  with $\eta_2:=\frac{e^{-2D}}{\abs{\mcw}}$. %\frac{C'}{\abs{\mcw}}e^{-2D}$.
  We now let $\underline{v}$ be the solution of
  \syss{
  \dr_t\underline{v}(t,w')&=(\Lambda-d(w)) \underline{v},&t>2\vg w'\in\mcw,\\
  \nu\cdot\mca\nabla \underline{v}(t,w')&=0,&w'\in\dr\mcw,\\
  \underline{v}(2,w')&\equiv\eta_2%\frac{\eta_2}{\abs{\mcw}}
  &w'\in\mcw.}
For all $w\in\mcw$ and for all $t\goq 2$,
  \begin{align*}
    \psf_{\delta_w}(\tau>t)
    &=\int_{\mcw}u^{(w)}(t,\cdot)
    \goq\int_{\mcw}\underline{v}(t,\cdot)
    = \eta_2\psf_{\mu}(\tau>t-2)
    \goq\eta_2\psf_{\mu}(\tau>t).
  \end{align*}
  Finally, for all $w\in\mcw$ and for all $t\in\cro{0,2}$, we have $\psf_{\delta_w}(\tau>t)\goq e^{-2D}\goq e^{-2D}\psf_{\mu}(\tau>t)$.
  Thus Assumption A2 holds with
  \[c_2:=\min\cro{\eta_2,e^{-2D}}>0.\]
  Hence, Assumptions A1 and A2 are satisfied with $\mu=\frac{\de w}{\abs{\mcw}}$. By Theorem 1.1 of~\cite{ChaVil14}, the first part of the theorem is shown. By Theorem 2.1 of~\cite{ChaVil14}, the constants $C$ and $\chi$ only depend on $t_0$, $c_1$ and $c_2$.
  We took $t_0=2$, and $c_1$ and $c_2$ only depend on $\underline{\mca}$, $\overline{\mca}$, ${\mcb}_M$, $\dabs{d}_{\infty}$ and $\mcw$. 
  Thus the constants $C$ and $\chi$ only depend on $\underline{\mca}$, $\overline{\mca}$, $\mcb_M$, $\dabs{d}_{\infty}$ and $\mcw$. 
\end{proof}

\section{Behaviour of the principal eigenfunction (Theorem~\ref{thm:diffusion_zero})}\label{s:proof_thm}

Throughout Section~\ref{s:proof_thm}, we will abuse notation and denote, for a measurable function $f:\mcw\to\Er$ and for measurable sets $E\subset\mcy$ and $F\subset\mcz$,
\[f(E,F)=\int_{E\times F}f(y,z)\de y\de z,\esp f(\cdot,F)=\int_{F}f(\cdot,z)\de z.\]
Throughout this section, we fix $y_0\in\overline{\mcy}$.

\subsection{Description of the framework}

Recall that in Subsection~\ref{ss:definition_process}, we defined a killed process $(Y^{\eps}_t,Z^{\eps}_t)_{t\goq 0}$ on $\mcw\cup\acc{\dr}$.
By Proposition~\ref{ppn:qsd}, we have, for all measurable sets $E\subset\mcy$ and $F\subset\mcz$,  for all $t>0$,
\begin{align*}
 \varphi_{\eps}(E,F)=\psf_{\varphi_{\eps}}\pth{({Y}^{\eps}_t,{Z}^{\eps}_t)\in E\times F\kt {\tau}^{\eps}>t}.
\end{align*}
With this in view, we understand why Theorem~\ref{thm:diffusion_zero} will be a consequence of the following proposition, the proof of which is the main part of Section~\ref{s:proof_thm}. For $\eta>0$, we set $V^{\eta}:=B(y_0,\eta)\cap\mcy$, where $B(y_0,\eta)\subset\Er^n$ is the ball with centre $y_0$ and radius $\eta$. 

\begin{ppn}\label{ppn:estimation_diff_zero}
  Let $y_0\in\overline{\mcy}$. For $\eps\in(0,1)$, set $t_{\eps}:=\frac{\ln(1/\eps)}{1+4(\max d-\min d)}$
  and $\eta_{\eps}:=\eps^{1/4}$. Let $E^{\eps}:=B\pth{y_0,\eps^2}\cap\mcy\subset V^{\eta_{\eps}}$.
  Then, as $\eps\to 0$,
  \begin{equation*}
    \dabs{\psf_{\varphi_\eps}\pth{{Z}^{\eps}_{t_{\eps}}\in  \de z\kt \tau^{\eps}>t_{\eps}\vg Y^{\eps}_{t_\eps}\in E^{\eps}}-\psi^{y_0}(z)\de z}_{TV}\to 0.
  \end{equation*}
  The convergence holds uniformly in $y_0\in\overline{\mcy}$.
\end{ppn}

Our goal is to focus on the behaviour of the process when $Y_t^{\eps}$ is close to $y_0$.
Therefore, for $t>0$, $\eta>0$, $\eps>0$, we let
\begin{align*}
  T(\eta,\eps;t)&:=\left\{
  \begin{aligned}
    &\inf\acc{s\in\ifo{0,t}\kt {Y}^{\eps}_{t-s}\notin V^{\eta}}&\text{if $\acc{s\in\ifo{0,t}\kt {Y}^{\eps}_{t-s}\notin V^{\eta}}\neq\emptyset$},\\
    &t&\text{otherwise,}
  \end{aligned}
  \right.
\end{align*}
be the time elapsed since $Y^{\eps}_t$ last entered the neighbourhood $V^{\eta}$ of $y_0$ without leaving it. See Figure~\ref{fig:reduction_env_hom}. The following lemma is the first step towards the proof of Proposition~\ref{ppn:estimation_diff_zero}.

\begin{lem}\label{lem:estimation_diff_zero}
There exist $K>0$, $\chi>0$ and $\chi'>0$, independent of $y_0$, such that for all $E\subset V^{\eta}$ and for all $t>0$, $\eta>0$, $\eps>0$,
\begin{align}
  &\dabs{\psf_{\varphi_{\eps}}\pth{{Z}^{\eps}_t\in \de z\kt{\tau}^{\eps}>t\vg Y_t^{\eps}\in E}-\psi^{y_0}(z)\de z}_{TV}\nonumber\\
  %&\dabs{\psf_{\varphi_{\eps}}\pth{{Z}^{\eps}_t\in \cdot\kt{\tau}^{\eps}>t\vg Y_t^{\eps}\in E}-\psi^{y_0}}_{TV}\nonumber\\
  &\esp\esp<K\esf_{\varphi_{\eps}}\cro{e^{-\chi T(\eta,\eps;t)}\kt \tau^{\eps}>t\vg Y^{\eps}_t\in E}+K\pth{e^{\eta\chi' t}-1}.\label{eq:estimation_expectations}
\end{align}
\end{lem}

\fig{10}{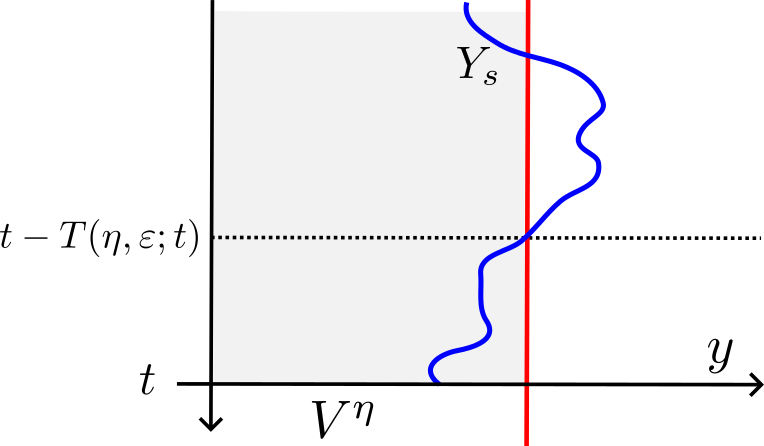}{\label{fig:reduction_env_hom} Definition of $T(\eta,\eps;t)$.}

The proof of Lemma~\ref{lem:estimation_diff_zero} is based on the observation that if we keep $Y^{\eps}_t$ near $y_0$, the quasi-stationary distribution for $Z^{\eps}_t$ is approximately $\psi^{y_0}$ (see Proposition~\ref{ppn:qsd}). The first term of the right-hand side of ~\eqref{eq:estimation_expectations} corresponds to how well the distribution of $Z^{\eps}_t$ conditionally on survival at time $t$ approximates the quasi-stationary distribution  $\psi^{y_0}$. We will see in Subsection~\ref{ss:frozen} how this term arises. The second term of the right-hand side of ~\eqref{eq:estimation_expectations} corresponds to the approximation made when forgetting the role of $Y^{\eps}_t$. We will see in Subsection~\ref{ss:comparison_frozen} how this term arises. Finally, in Subsection~\ref{ss:conclusion_big_thm}, we will conclude the proof of Proposition~\ref{ppn:estimation_diff_zero} and the proof of Theorem~\ref{thm:diffusion_zero}.

% 
%\fig{7}{img/def_processus_qsd.png}{\label{fig:def_processus_qsd}\wtd need update The process $({Y}_s)_s$ is a drifted Brownian motion in $\Er$ and corresponds to the position of a particle. The process $({Z}_s)_s$ is a Brownian motion in the bounded set $\mcz$ and corresponds to the state of the particle. The trajectory is from left to right.}

\subsection{Work with a locally frozen death rate}\label{ss:frozen}

In Subsections~\ref{ss:frozen} and~\ref{ss:comparison_frozen}, we will work with a fixed $\eps>0$. In order to lighten the notations, we will therefore omit the superindex $\eps$. The notations become:
\[\widetilde Y=\widetilde Y^{\eps},\esp \widetilde Z=\widetilde Z^{\eps},\esp \tau=\tau^{\eps},\esp T(\eta;t)=T(\eta,\eps;t),\esp \varphi=\varphi_{\eps}, \qquad\text{etc}.\]
We should not forget, however, that all these objects do indeed depend on $\eps$.

We call $(\widetilde{Y}_t,\widetilde{Z}_t)_{t\goq 0}$ the carrying process; in the sequel, we will construct several processes based on $(\widetilde{Y}_t,\widetilde{Z}_t)_{t\goq 0}$ but dying at different stopping times.

\paragraph{Approximate the stopping time $\tau$ by another stopping time $\tau^{\eta}$.}

Let $\eta>0$ and let $V^{\eta}:=B(y_0,\eta)\cap\mcy$. We consider a function $d^{\eta}$ which approximates $d$ and which, on the neighbouring band $V^{\eta}\times\mcz$, only depends on $z$:
\begin{align*}\label{eq:dr_homogeneous_near_y0}
  d^{\eta}(y,z)&=d(y,z),& y\notin V^{\eta},\\
  d^{\eta}(y,z)&=d(y_0,z),& y\in V^{\eta}.
\end{align*}
See Figure~\ref{fig:constant_near_y0}.

We define a stochastic process $(Y^{\eta}_t,Z^{\eta}_t)_{t\goq 0}$ exactly as $(Y_t,Z_t)_{t\goq 0}$, but which dies at rate $d^{\eta}$ instead of dying at rate $d$. Namely, we let $\tau^{\eta}$ be a clock with rate $(d^{\eta}(\widetilde{Y}_t,\widetilde{Z}_t))_{t\goq 0}$ and we set
  \[
  ({Y}^{\eta}_t,{Z}^{\eta}_t):=
  \left\{
  \begin{aligned}
    &(\widetilde{Y}_t,\widetilde{Z}_t),&t<{\tau}^{\eta},\\
    &\dr,&t\goq{\tau}^{\eta}.
  \end{aligned}
  \right.
  \]
\fig{8}{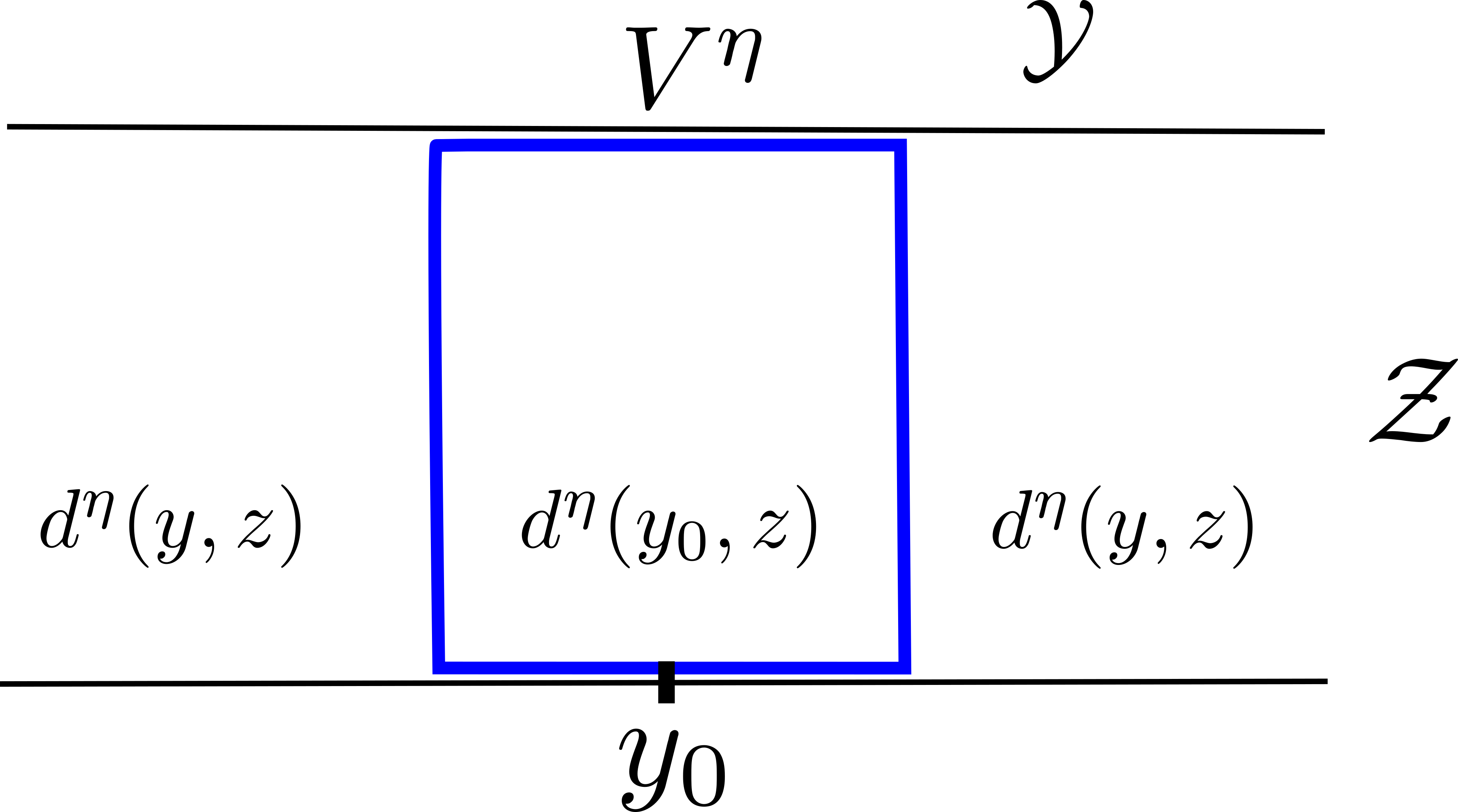}{\label{fig:constant_near_y0}In the neighbourhood $V^{\eta}$ of $y_0$, $d^{\eta}$ is independent of $y$, but may still depend on~$z$.}

\paragraph{Introduce processes $(\widehat{Y}_s,\widehat{Z}_s)_{s\goq 0}$ and $(\widehat{Y}^{\eta}_s,\widehat{Z}^{\eta}_s)_{s\goq 0}$ living in $V^{\eta}\times\mcz\subset\mcw$, and their stopping times $\widehat{\tau}$ and $\widehat{\tau}^{\eta}$.}

We can bear in mind that between $t-T(\eta;t)$ and $t$, the coefficients of the SDE analogous to~\eqref{eq:sde} satisfied by $Z^{\eta}$ are independent of $Y^{\eta}$. To formalise this idea, we consider auxiliary processes $(\widehat{Y}_s,\widehat{Z}_s)_s$ and $(\widehat{Y}^{\eta}_s,\widehat{Z}^{\eta}_s)_s$ which behave in a similar way as $(Y_t,Z_t)_t$ and $({Y}^{\eta}_t,{Z}^{\eta}_t)_t$, but which die upon reaching $\dr V^{\eta}\times\mcz$.

Namely, we let $\widehat{\tau}_{d}$ and $\widehat{\tau}_d^{\eta}$ be two clocks with respective rates $d(\widetilde{Y}_s,\widetilde{Z}_s)$ and $d^{\eta}(\widetilde{Y}_s,\widetilde{Z}_s)$. 
We let also:
\begin{align*}
  \widehat{\tau}_Y&:=\inf\acc{u \goq 0\tq \widetilde{Y}_u\notin V^{\eta}}
\end{align*}
  be the first exit time from $V^{\eta}$ by $\widetilde{Y}$. 
  Finally, we set $\widehat{\tau}:=\widehat{\tau}_{d}\wedge \widehat{\tau}_Y$ and $\widehat{\tau}^{\eta}:=\widehat{\tau}_{d}^{\eta}\wedge \widehat{\tau}_Y$.
  We now define
  \[
  (\widehat{Y}_s,\widehat{Z}_s):=
  \left\{
  \begin{aligned}
    &(\widetilde{Y}_s,\widetilde{Z}_s),&t<\widehat{\tau},\\
    &\dr,&t\goq\widehat{\tau};
  \end{aligned}
  \right.
  \esp\esp(\widehat{Y}^{\eta}_s,\widehat{Z}^{\eta}_s):=
  \left\{
  \begin{aligned}
    &(\widetilde{Y}^{\eta}_s,\widetilde{Z}^{\eta}_s),&t<\widehat{\tau}^{\eta},\\
    &\dr,&t\goq\widehat{\tau}^{\eta}.
  \end{aligned}
  \right.
  \]
  Hence $\widehat{\tau}$ and $\widehat{\tau}^{\eta}$ are the respective death times of the processes $(\widehat{Y}_u,\widehat{Z}_u)_u$ and $(\widehat{Y}_u^{\eta},\widehat{Z}_u^{\eta})_u$.
  
  We will first compare $(\widehat{Y}_s,\widehat{Z}_s)_s$ and $({Y}_t,{Z}_t)_t$.
    Then, we will see how to approximate the local principal eigenfunctions
    $\psi^{y}$, $y\in\overline{\mcy}$,
    using the process $(\widehat{Y}_s^{\eta},\widehat{Z}_s^{\eta})_{s}$.
    For this approximation to hold,
    we need to use the independence of the processes
    $(\widehat{Y}^{\eta}_s)_s$ and $(\widehat{Z}^{\eta}_s)_s$
    when started from some $y\in V^{\eta}$ and $z\in\mcz$ respectively;
    to get the independence, we will use Assumption~\eqref{eq:dependence_restriction}.
    Finally, in the next subsection,
    we will compare the processes $(\widehat{Y}_s,\widehat{Z}_s)_s$
    and $(\widehat{Y}_s^{\eta},\widehat{Z}_s^{\eta})_{s}$,
    and prove Lemma~\ref{lem:estimation_diff_zero}.

\medskip

For integer $k$ with $0\loq k<t$, we set $I_k:=\ifo{k,k+1}$.
We let $m_k(t)$ be the law of $({Y}_{t-k}\vg {Z}_{t-k})$ conditionally on the event $\acc{T(\eta;t)\in I_k,\tau>t}$. That is, for $E\subset\mcy$ and $F\subset\mcz$,
\[m_k(t)\pth{E\times F}:=\psf_{\varphi}\pth{({Y}_{t-k}\vg {Z}_{t-k})\in E\times F\kt T(\eta;t)\in I_k\vg \tau>t}.\]
The following lemma compares the two processes $(\widehat{Y}_s,\widehat{Z}_s)_s$ and $({Y}_t,{Z}_t)_t$.

\begin{lem}\label{lem:equality_y_yhat}
  For all measurable subsets $E\subset V^{\eta}$ and $F\subset\mcz$, for all $t>0$ and integer $k$ with $0\loq k<t$,
  \begin{equation*}\label{eq:equality_y_yhat}
  \psf_{\varphi}\pth{({Y}_t,{Z}_t)\in E\times F\kt \tau>t\vg T(\eta;t)\in I_k}=\psf_{m_k(t)}\pth{(\widehat{Y}_k,\widehat{Z}_k)\in E\times F\kt\widehat{\tau}>k}.
\end{equation*}  
\end{lem}

\begin{proof}
  Let $t_0:=t-k$. We have:
  \begin{align*}
    &\psf_{\varphi}\pth{({Y}_t,{Z}_t)\in E\times F\kt \tau>t\vg T(\eta;t)\in I_k}\\
    &\esp=\psf_{\varphi}\pth{(\widetilde{Y}_{t_0+k},\widetilde{Z}_{t_0+k})\in E\times F\kt \tau>t_0+k\vg T(\eta;t)\in I_k}.
  \end{align*}
  Now, we write:
    \begin{align*}
      &\psf_{\varphi}\pth{(\widetilde{Y}_{t_0+k},\widetilde{Z}_{t_0+k})\in E\times F\kt \tau>t_0+k\vg T(\eta;t)\in I_k}\\
      &\qquad=\int_{\mcy\times\mcz}\psf_{\varphi}\pth{(\widetilde{Y}_{t_0},\widetilde{Z}_{t_0})\in \de y\de z,\, (\widetilde{Y}_{t_0+k},\widetilde{Z}_{t_0+k})\in E\times F\kt \tau>t_0+k\vg T(\eta;t)\in I_k}\\
      &\qquad=\int_{\mcy\times\mcz}m_k(t)(\de y\de z)\psf_{\varphi}\pth{(\widetilde{Y}_{t_0+k},\widetilde{Z}_{t_0+k})\in E\times F\kt (\widetilde{Y}_{t_0},\widetilde{Z}_{t_0})\in \de y\de z,\,\tau>t_0+k\vg T(\eta;t)\in I_k},
    \end{align*}
    where we note
    \begin{multline*}
      \psf_{\varphi}\pth{(\widetilde{Y}_{t_0+k},\widetilde{Z}_{t_0+k})\in \cdot\kt (\widetilde{Y}_{t_0},\widetilde{Z}_{t_0})\in \de y\de z,\,\tau>t_0+k\vg T(\eta;t)\in I_k}\\
      :=\lim_{\eps\to0}\psf_{\varphi}\pth{(\widetilde{Y}_{t_0+k},\widetilde{Z}_{t_0+k})\in \cdot\kt (\widetilde{Y}_{t_0},\widetilde{Z}_{t_0})\in B((y,z),\eps),\,\tau>t_0+k\vg T(\eta;t)\in I_k}.
    \end{multline*}
The validity of this expression and of the computations follows from the fact that $(\widetilde{Y}_{t},\widetilde{Z}_{t})$ has a continuous density for all $t\goq0$.

Next, let us rewrite the following events:
    \[\acc{T(\eta;t)\in I_k}=\acc{\exists u\in[t_0-1,t_0],\,\widetilde{Y}_{u}\notin V^{\eta}}\cap\acc{\forall u\in[t_0,t],\, \widetilde{Y}_{u}\in V^{\eta}}=:A_1\cap A_2,\]
    and
    \[\acc{\tau>t_0+k}=\acc{\tau>t_0}\cap \acc{\tau'>k}=:B_1\cap B_2,\]
    where~$\tau'$ is a clock with rate $d(\widetilde{Y}_{t_0+t},\widetilde{Z}_{t_0+t})_{t\goq 0}$ (which satisfies $\tau'=\tau-k$ on $\acc{\tau>k}$).
    Now, we use the Markov property, which implies that conditionnally on $(\widetilde{Y}_{t_0},\widetilde{Z}_{t_0})$, $(\widetilde{Y}_{t_0+k},\widetilde{Z}_{t_0+k})$ is independent of $(\widetilde{Y}_{u},\widetilde{Z}_{u})_{u\in[0,t_0]}$, and therefore is independent of $A_1$ and $B_1$.
    We obtain:
    \begin{multline*}
      \psf_{\varphi}\pth{(\widetilde{Y}_{t_0+k},\widetilde{Z}_{t_0+k})\in E\times F\kt (\widetilde{Y}_{t_0},\widetilde{Z}_{t_0})\in \de y\de z,\,\tau>t_0+k\vg T(\eta;t)\in I_k}\\
      =\psf_{(y,z)}\pth{(\widetilde{Y}'_{k},\widetilde{Z}'_{k})\in E\times F\kt A_2,\, B_2},
    \end{multline*}
    where $(\widetilde{Y}'_{t},\widetilde{Z}'_{t}):=(\widetilde{Y}'_{t_0+t},\widetilde{Z}'_{t_0+t})$.
    Gathering everything together, we obtain:
    \begin{align*}
    &\psf_{\varphi}\pth{({Y}_t,{Z}_t)\in E\times F\kt \tau>t\vg T(\eta;t)\in I_k}\\
      &\esp=\psf_{m_k(t)}\pth{(\widetilde{Y}'_k,\widehat{Z}'_k)\in E\times F\kt A_2,B_2}\\  
      &\esp=\psf_{m_k(t)}\pth{(\widetilde{Y}_k,\widetilde{Z}_k)\in E\times F\kt \widehat{\tau}_d>k\vg \widehat{\tau}^Y>k}\\
    &\esp=\psf_{m_k(t)}\pth{(\widehat{Y}_k,\widehat{Z}_k)\in E\times F\kt \widehat{\tau}>k}.
  \end{align*}
\end{proof}

\paragraph{Introduce a process $(\overline{Z}_s)_{s\goq 0}$ living in $\mcz$.}
We consider another auxiliary process $(\overline{Z}_s)_s$ living in $\mcz$, which behaves in a similar way as $({Z}_t)_{t\goq 0}$ and $(\widehat{Z}_s)_{s\goq 0}$, but independently of any process in~$\mcy$.
We let $\overline{\tau}$ be a clock with rate $(d(y_0,\widetilde{Z_s}))_{s\goq 0}$. We define
\[\overline{Z}_t:=
\left\{
\begin{aligned}
  &\widetilde{Z}_t&t<\overline{\tau},\\
  &\dr&t\goq \overline{\tau}.
\end{aligned}
\right.
\]
Let us try to understand the behaviour of $(\overline{Z}_s)_{s\goq 0}$. In~\eqref{eq:pde_evol_yz}, if we take the initial condition~$m$ independent of $y$, then the solution $\widetilde{u}^{\eps}$ remains independent of $y$ (due to~\eqref{eq:dependence_restriction}). In this case, therefore, the density of $\widetilde{Z}_t$ is described by the solution $\widetilde{v}$ of the equation
\syss{
  \dr_t\widetilde{v}(t,z)&=\mcl_z\widetilde{v},&t\goq 0\vg z\in\mcz,\\
  \nu\cdot a\nabla\widetilde{v}&=0,&z\in\dr\mcz,\\
  \widetilde{v}(0,z)&=m(z),&z\in\mcz.
}
Taking $\mcw=\mcz$ and $\Lambda=\mcl_z$, we retrieve the setting of Section~\ref{s:qsd}. By Proposition~\ref{ppn:qsd}, then, we find that the quasi-stationary distribution of the process $(\overline{Z}_s)_{s\goq 0}$ is $\psi^{y_0}$. 
Also, note that $(\overline{Z}_s)_s$ and $(\widehat{Z}^{\eta}_s)_s$ have much the same behaviour. Therefore, in view of Theorem~\ref{thm:cv_qsd}, the following lemma is natural.

\begin{lem}\label{lem:cv_shat}
  There exist $C>0$ and $\chi>0$ such that for all initial distribution $m$ on $\mcw$ and for all $E\subset V^{\eta}$,
  \[\dabs{\psf_{m}\pth{\widehat{Z}^{\eta}_s\in\de z\kt\widehat{\tau}^{\eta}>s\vg\widehat{Y}^{\eta}_s\in E}-\psi^{y_0}(z)\de z}_{TV}\loq Ce^{-\chi s}.\]
  The constants $C>0$ and $\chi>0$ can be chosen independently of $y_0$.
\end{lem}

\begin{proof}
  %The result is a consequence of~\eqref{eq:comp_shat_sbar} and of Theorem~\ref{thm:cv_qsd} applied to the process $(\overline{Z}_u)_u$. Let us prove~\eqref{eq:comp_shat_sbar}.

  We have, for any initial distribution $m$ on $\mcw$,
  \begin{align*}
    \psf_{m}\pth{\widehat{Y}^{\eta}_s\in E,\widehat{Z}^{\eta}_s\in\cdot\kt \widehat{\tau}^{\eta}>s}&=\frac{\psf_{m}\pth{\widehat{Y}^{\eta}_s\in E,\widehat{Z}^{\eta}_s\in\cdot\vg\widehat{\tau}^{\eta}>s}}{\psf_{m}\pth{\widehat{\tau}^{\eta}>s}}\\
    &=\frac{1}{\psf_{m}\pth{\widehat{\tau}^{\eta}>s}}\int_{\mcw}\psf_{(y,z)}\pth{\widehat{Y}^{\eta}_s\in E,\widehat{Z}^{\eta}_s\in\cdot\vg\widehat{\tau}^{\eta}>s}m(\de y,\de z).
  \end{align*}
  For conciseness, we denoted $\psf_{(y,z)}$ instead of $\psf_{\delta_{(y,z)}}$.
  Let $(y,z)\in V^{\eta}$. The properties of the Cauchy problem~\eqref{eq:cp_evol_wtilde} imply that the processes $(\widetilde{Y}_s)_{s\goq 0}$ starting from $y$ and $(\widetilde{Z}_s)_{s\goq 0}$ starting from $z$ are independent (this is where we need the restrictions~\eqref{eq:dependence_restriction} on the dependencies of the coefficients).
  Moreover, the stopping time~$\widehat{\tau}^{\eta}_Y$ is a function of the process $(\widetilde{Y}_s)_{s\goq 0}$ only, and the stopping time $\widehat{\tau}^{\eta}_d$ is a function of the process $(\widetilde{Z}^{\eta}_s)_{s\goq 0}$ only. Therefore, for all $(y,z)\in\mcw$, for all $s>0$,
  \begin{align*}
    \psf_{(y,z)}\pth{\widehat{Y}^{\eta}_s\in E,\widehat{Z}^{\eta}_s\in\cdot\vg \widehat{\tau}^{\eta}>s}=\psf_{(y,z)}\pth{\widehat{Z}^{\eta}_s\in\cdot\vg \widehat{\tau}^{\eta}_{d}>s}\psf_{(y,z)}\pth{\widehat{Y}^{\eta}_s\in E,\widehat{\tau}^{\eta}_{Y}>s}.
  \end{align*}
  We obtain:
  \begin{align*}
    &\psf_{m}\pth{\widehat{Y}^{\eta}_s\in E,\widehat{Z}^{\eta}_s\in\cdot\kt \widehat{\tau}^{\eta}>s}\\
    &\esp=\frac{1}{\psf_{m}\pth{\widehat{\tau}^{\eta}>s}}\int_{\mcw}\psf_{(y,z)}\pth{\widehat{Z}^{\eta}_s\in\cdot\vg\widehat{\tau}^{\eta}_{d}>s}\psf_{(y,z)}\pth{\widehat{Y}^{\eta}_s\in E,\widehat{\tau}^{\eta}_{Y}>s}m(\de y,\de z)\\
    &\esp=\frac{1}{\psf_{m}\pth{\widehat{\tau}^{\eta}>s}}\int_{\mcw}\psf_{(y,z)}\pth{\widehat{Z}^{\eta}_s\in\cdot\kt\widehat{\tau}^{\eta}_{d}>s}\psf_{(y,z)}\pth{\widehat{\tau}^{\eta}_{d}>s}\psf_{(y,z)}\pth{\widehat{Y}^{\eta}_s\in E,\widehat{\tau}^{\eta}_{Y}>s}m(\de y,\de z)\\
    &\esp=\frac{1}{\psf_{m}\pth{\widehat{\tau}^{\eta}>s}}\int_{\mcw}\psf_{(y,z)}\pth{\widehat{Z}^{\eta}_s\in\cdot\kt\widehat{\tau}^{\eta}_{d}>s}\psf_{(y,z)}\pth{\widehat{Y}^{\eta}_s\in E,\widehat{\tau}^{\eta}>s}m(\de y,\de z).
  \end{align*}
  %using the fact that conditionally on $\acc{Y_0=y,Z_0=z}$, the stopping times $\widehat{\tau}^{\eta}_{d}$ is independent of $(\widetilde{Y}_t)_{t\goq 0}$ and $\widehat{\tau}^{\eta}_{Y}$.
  Therefore,
  \begin{align*}
    &\dabs{\psf_{m}\pth{\widehat{Y}^{\eta}_s\in E,\widehat{Z}^{\eta}_s\in \de z\kt \widehat{\tau}^{\eta}>s}-\psf_{m}\pth{\widehat{Y}^{\eta}_s\in E\kt\widehat{\tau}^{\eta}>s}\psi^{y_0}(z)\de z}_{TV}\\
    &\esp\loq \frac{1}{\psf_{m}\pth{\widehat{\tau}^{\eta}>s}}\int_{\mcw}\dabs{\psf_{(y,z)}\pth{\widehat{Z}^{\eta}_s\in\de z'\kt\widehat{\tau}^{\eta}_{d}>s}-\psi^{y_0}(z')\de z'}_{TV}\psf_{(y,z)}\pth{\widehat{Y}^{\eta}_s\in E,\widehat{\tau}^{\eta}>s}m(\de y,\de z).
  \end{align*}
  But, by Theorem~\ref{thm:cv_qsd} applied to the killed process $(\overline{Z}_s)_{s\goq 0}$ (see also the discussion before the statement of Lemma~\ref{lem:cv_shat}),
  \begin{align*}
    \dabs{\psf_{(y,z)}\pth{\widehat{Z}^{\eta}_s\in\de z\kt\widehat{\tau}^{\eta}_{d}>s}-\psi^{y_0}(z)\de z}_{TV}&=\dabs{\psf_{z}\pth{\overline{Z}_s\in\de z\kt\overline{\tau}>s}-\psi^{y_0}(z)\de z}_{TV}\\
    &\loq Ce^{-\chi s}.
  \end{align*}
  Note that the constants $C$ and $\chi$ can be chosen independently of $y_0$, by the second statement in Theorem~\ref{thm:cv_qsd}.
  Hence
  \begin{align*}
    &\dabs{\psf_{m}\pth{\widehat{Y}^{\eta}_s\in E,\widehat{Z}^{\eta}_s\in \de z\kt \widehat{\tau}^{\eta}>s}-\psf_{m}\pth{\widehat{Y}^{\eta}_s\in E\kt\widehat{\tau}^{\eta}>s}\psi^{y_0}(z)\de z}_{TV}\\
    &\esp\loq Ce^{-\chi s} \times\frac{1}{\psf_{m}\pth{\widehat{\tau}^{\eta}>s}}\int_{\mcw}\psf_{(y,z)}\pth{\widehat{Y}^{\eta}_s\in E,\widehat{\tau}^{\eta}>s}m(\de y,\de z)\\
    &\esp= Ce^{-\chi s}\psf_{m}\pth{\widehat{Y}^{\eta}_s\in E\kt\widehat{\tau}^{\eta}>s}.
  \end{align*}
  Dividing both sides by $\psf_{m}\pth{\widehat{Y}^{\eta}_s\in E\kt\widehat{\tau}^{\eta}>s}$, we obtain:
  \[\dabs{\psf_{m}\pth{\widehat{Z}^{\eta}_s\in \de z\kt \widehat{\tau}^{\eta}>s\vg \widehat{Y}^{\eta}_s\in E}-\psi^{y_0}(z)\de z}_{TV}\loq Ce^{-\chi s}.\]
Finally, we mentioned that the constants $C$ and $\chi$ can be chosen independently of $y_0$.
\end{proof}

\subsection{Comparison between $\widehat{\tau}$ and $\widehat{\tau}^{\eta}$}\label{ss:comparison_frozen}

The processes $(\widehat{Y}_t,\widehat{Z}_t)_{t\goq 0}$ and $(\widehat{Y}_t^{\eta},\widehat{Z}_t^{\eta})_{t\goq 0}$, and their respective death times $\widehat{\tau}$ and $\widehat{\tau}^{\eta}$, are almost the same. The only difference is that the dependence in $y$ of the death rate is removed in the definition of $(\widehat{Y}^{\eta}_s,\widehat{Z}^{\eta}_s)$. With this respect, $(\widehat{Y}^{\eta}_s,\widehat{Z}^{\eta}_s)$ dies at a rate \enquote{frozen} at $y_0$ while $(\widehat{Y}_s,\widehat{Z}_s)$ dies at a rate depending on the $y$-coordinate.
For small $\eta$, the coefficient $d$ does not depend much on $y\in V^{\eta}$, in the sense that for all $y\in V^{\eta}$,
\[d(y,\cdot)\simeq d(y_0,\cdot).\]
Therefore, we expect that $\widehat{\tau}\simeq\widehat{\tau}^{\eta}$.
This intuition is made precise in the following lemma, which is the last step before the proof of Lemma~\ref{lem:estimation_diff_zero}.

\begin{lem}\label{lem:comparison_frozen}
There exists $\chi'>0$ such that for all $s>0$, for all $E\subset V^{\eta}$, for all probability distribution $m$ on $V^{\eta}$,
\[    \dabs{\psf_{m}\pth{{\widehat{Z}}_s\in\cdot\kt\widehat{\tau}>s\vg \widehat{Y}_s\in E}-\psf_{m}\pth{{\widehat{Z}}^{\eta}_s\in\cdot\kt\widehat{\tau}^{\eta}>s\vg \widehat{Y}^{\eta}_s\in E}}_{TV}
\loq 2\pth{e^{\eta\chi' t}-1}.\]
The constant $\chi'$ can be chosen independently of $y_0$.
\end{lem}

\begin{proof}
  Recall that $\widehat{\tau}=\widehat{\tau}_{d}\wedge \widehat{\tau}_Y$ where $\widehat{\tau}_{d}$ depends on the death rate $d(\widetilde{Y}_t,\widetilde{Z}_t)$, and $\widehat{\tau}_Y$ is the time at which $\widetilde{Y}$ reaches the boundary $\dr V^{\eta}$. Likewise, $\widehat{\tau}^{\eta}=\widehat{\tau}^{\eta}_{d}\wedge \widehat{\tau}_Y$. We therefore decompose:
  \begin{align*}
    \acc{{\widehat{\tau}}>s}=\acc{{\widehat{\tau}}_{d}>s\vg{\widehat{\tau}}_Y>s },\esp    \acc{{\widehat{\tau}}^{\eta}>s}=\acc{{\widehat{\tau}}^{\eta}_{d}>s\vg{\widehat{\tau}}_Y>s }.
  \end{align*}
  Let $\kappa>0$ be the Lipschitz constant of $d$ in the $y$-direction:
  \[\kappa:=\sup_{z\in\mcz}\pth{\sup_{y,y'\in\mcy\vg y\neq y'}\frac{\abs{d(y,z)-d(y',z)}}{\abs{y-y'}}}<+\infty.\]
  On the event $\acc{{\widehat{\tau}}_Y>s}$, we then have, for all $u\in\cro{0,s}$,
  \[\abs{d({\widehat{Y}}_u,{\widehat{Z}}_u)-d(y_0,{\widehat{Z}}_u)}\loq \kappa\eta.\]
  For conciseness of notations, denote by $\mcd$ the event $\mcd=\acc{\widehat{\tau}_Y>s\vg \widehat{Y}_s\in E}$. Recall the definition of $\widehat{\tau}_d$ in Subsection~\ref{ss:frozen}, using Definition~\ref{dfi:poisson}. %For $\lambda>0$, let $Poi(\lambda)$ denote a Poisson random variable with parameter $\lambda$.
  We have:
  \begin{align*}
    \psf_{m}\pth{{\widehat{\tau}}_{d}>s\kt \mcd}
    &=\esf_{m}\cro{\exp\pth{-\int_0^sd({\widetilde{Y}}_u,{\widetilde{Z}}_u)\de u}\kt \mcd}\\
    &\goq\esf_{m}\cro{\exp\pth{-\kappa\eta s-\int_0^sd(y_0,{\widetilde{Z}}_u)\de u}\kt \mcd}.
  \end{align*}
  This gives:
  \begin{align*}
    \psf_{m}\pth{{\widehat{\tau}}_{d}>s\kt \mcd}
    &\goq e^{-\kappa\eta s}\esf_{m}\cro{\exp\pth{-\int_0^sd(y_0,{\widetilde{Z}}_u)\de u}\kt \mcd}\\
    &=e^{-\kappa\eta s}\psf_{m}\pth{\widehat{\tau}^{\eta}_{d}>s\kt \mcd}.
  \end{align*}
%  We have:
%  \begin{align*}
%    \psf_{m}\pth{{\widehat{\tau}}_{d}>s\kt \mcd}
%    &=\esf_{m}\cro{\prb{Poi\pth{\int_0^sd({\widehat{Y}}_u,{\widehat{Z}}_u)\de u}=0}\kt \mcd}\\
%    &\goq\esf_{m}\cro{\prb{Poi\pth{\kappa\eta s+\int_0^sd(y_0,{\widehat{Z}}_u)\de u}=0}\kt \mcd}.
%  \end{align*}
%  This gives:
%  \begin{align*}
%    \psf_{m}\pth{{\widehat{\tau}}_{d}>s\kt \mcd}
%    &\goq e^{-\kappa\eta s}\esf_{m}\cro{\prb{Poi\pth{\int_0^sd(y_0,{\widehat{Z}}_u)\de u}=0}\kt \mcd}\\
%    &=e^{-\kappa\eta s}\psf_{m}\pth{\widehat{\tau}^{\eta}_{d}>s\kt \mcd}.
%  \end{align*}
  Likewise, $\psf_{m}\pth{\widehat{\tau}^{\eta}_{d}>s\kt \mcd}\goq e^{-\kappa\eta s}\psf_{m}\pth{\widehat{\tau}_{d}>s\kt \mcd}$. Thus:
  \[e^{-\kappa\eta s}\psf_{m}\pth{\widehat{\tau}^{\eta}_{d}>s\kt \mcd}\loq \psf_{m}\pth{\widehat{\tau}_{d}>s\kt \mcd}\loq e^{\kappa\eta s}\psf_{m}\pth{\widehat{\tau}^{\eta}_{d}>s\kt \mcd}.\]
  Similarly, for all $F\subset\mcz$, 
  \begin{align*}
    e^{-\kappa\eta s}\psf_{m}\pth{{\widehat{Z}}^{\eta}_s\in F,\widehat{\tau}^{\eta}_{d}>s\kt \mcd}
    &\loq \psf_{m}\pth{{\widehat{Z}}_s\in F,\widehat{\tau}_{d}>s\kt \mcd}\loq e^{\kappa\eta s}\psf_{m}\pth{{\widehat{Z}}^{\eta}_s\in F,\widehat{\tau}^{\eta}_{d}>s\kt \mcd}.    
  \end{align*}
  Therefore, for all $F\subset\mcz$,
  \begin{align*}
    e^{-2\kappa\eta s}\psf_{m}\pth{{\widehat{Z}}_s^{\eta}\in F\kt \widehat{\tau}^{\eta}_{d}>s\vg \mcd}
    &\loq \psf_{m}\pth{{\widehat{Z}}_s\in F\kt \widehat{\tau}_{d}>s\vg\mcd}\\
    &\loq e^{2\kappa\eta s}\psf_{m}\pth{{\widehat{Z}}_s^{\eta}\in F\kt \widehat{\tau}^{\eta}_{d}>s\vg \mcd}.
  \end{align*}
  Since this holds for all $F\subset\mcz$, we get:
   \begin{align*}
     &\dabs{\psf_{m}\pth{{\widehat{Z}}_s\in\cdot\kt \widehat{\tau}_{d}>s\vg\mcd}-\psf_{m}\pth{{\widehat{Z}}^{\eta}_s\in\cdot\kt  \widehat{\tau}^{\eta}_{d}>s\vg\mcd}}_{TV}\\
     &\esp\loq 2\,\max(e^{2\kappa\eta s}-1\vg 1-e^{-2\kappa\eta s})=2\pth{e^{2\kappa\eta s}-1}.
   \end{align*}
   % Consider $\tau_m=\widehat{\tau}\wedge\widehat{\tau}^{\eta}$ and let $(\widehat{Y}'_s)_{s\goq 0}$ and $(\widehat{Y}^{\eta}'_s)_{s\goq 0}$ be two processes equal respectively to the processes $(\widehat{Y}_s)_{s\goq 0}$ and $(\widehat{Y}^{\eta}_s)_{s\goq 0}$ for $s<\tau_m$, and be equal to $\dr$ for $s\goq \tau_m$.
   Finally, note that before their death, the processes $(\widehat{Y}_s,\widehat{Z}_s)_{s\goq 0}$ and $(\widehat{Y}^{\eta}_s,\widehat{Z}^{\eta}_s)_{s\goq 0}$ coincide with $(\widetilde{Y}_s,\widetilde{Z}_s)_{s\goq 0}$.
   In particular, $\acc{\widehat{\tau}^{\eta}_{d}>s,\mcd}=\acc{\widehat{\tau}^{\eta}_{d}>s,\widehat{\tau}_Y>s\vg \widetilde{Y}_s\in E}$.
   These elements yield
   \begin{align*}
     \psf_{m}\pth{{\widehat{Z}}^{\eta}_s\in\cdot\kt  \widehat{\tau}^{\eta}_{d}>s\vg\mcd}
     &=\psf_{m}\pth{{\widehat{Z}}^{\eta}_s\in\cdot\kt  \widehat{\tau}^{\eta}_{d}>s\vg\widehat{\tau}_Y>s\vg \widetilde{Y}_s\in E}\\
     &=\psf_{m}\pth{{\widehat{Z}}^{\eta}_s\in\cdot\kt  \widehat{\tau}^{\eta}_{d}>s\vg\widehat{\tau}_Y>s\vg\widehat{Y}^{\eta}_s\in E}\\
     &=\psf_{m}\pth{{\widehat{Z}}^{\eta}_s\in\cdot\kt  \widehat{\tau}^{\eta}>s\vg\widehat{Y}^{\eta}_s\in E}.
   \end{align*}
   We obtain: for all $s>0$,
   \begin{align*}
     &\dabs{\psf_{m}\pth{{\widehat{Z}}_s\in\cdot\kt\widehat{\tau}>s\vg\widehat{Y}_s\in E}-\psf_{m}\pth{{\widehat{Z}}^{\eta}_s\in\cdot\kt\widehat{\tau}^{\eta}>s\vg\widehat{Y}^{\eta}_s\in E}}_{TV}\\
     &\esp=\dabs{\psf_{m}\pth{{\widehat{Z}}_s\in\cdot\kt \widehat{\tau}_{d}>s\vg\mcd}-\psf_{m}\pth{{\widehat{Z}}_s\in\cdot\kt  \widehat{\tau}^{\eta}_{d}>s\vg\mcd}}_{TV}\\
     &\esp\loq e^{2\kappa\eta s}-1.
   \end{align*}
   Thus the lemma holds with $\chi'=2\kappa$, which is independent of $y_0$.
\end{proof}

We are ready for end of the proof of Lemma~\ref{lem:estimation_diff_zero}.

\begin{proof}[Proof of Lemma~\ref{lem:estimation_diff_zero}]
  Let us decompose the following probability according to the value of~$T(\eta;t)$: for all $t>0$,
  \begin{align*}
    &\psf_{\varphi}\pth{Y_t\in E,Z_t\in \cdot\kt \tau>t}\\
    &\esp=\sum_{k=0}^{\floor{t}}\psf_{\varphi}\pth{Y_t\in E,T(\eta;t)\in I_k\kt \tau>t}\psf_{\varphi}\pth{Z_t\in \cdot\kt \tau>t\vg Y_t\in E,T(\eta;t)\in I_k}\\
    &\esp=\sum_{k=0}^{\floor{t}}\psf_{\varphi}\pth{Y_t\in E,T(\eta;t)\in I_k\kt \tau>t}\psf_{m_k}\pth{\widehat{Z}_k\in \cdot\kt \widehat{\tau}>k\vg \widehat{Y}_k\in E},
  \end{align*}
  by Lemma~\ref{lem:equality_y_yhat} (using the shortcut $m_k=m_k(t)$).
  Thus, by the triangle inequality,
  \begin{align*}
    &\dabs{\psf_{\varphi}\pth{Y_t\in E,Z_t\in \de z\kt \tau>t}-\psf_{\varphi}(Y_t\in E\kt \tau>t)\psi^{y_0}(z)\de z}_{TV}\\
    &\esp\loq\sum_{k=0}^{\floor{t}}\psf_{\varphi}\pth{Y_t\in E,T(\eta;t)\in I_k\kt \tau>t}\dabs{\psf_{m_k}\pth{\widehat{Z}_k\in \de z\kt \widehat{\tau}>k\vg \widehat{Y}_k\in E}-\psi^{y_0}(z)\de z}_{TV}.
  \end{align*}
  By the triangle inequality again,
  and then Lemmas~\ref{lem:comparison_frozen} and~\ref{lem:cv_shat}, we have for all $k\goq 0$:
  \begin{align*}
    &\dabs{\psf_{m_k}\pth{\widehat{Z}_k\in \de z\kt \widehat{\tau}>k\vg \widehat{Y}_k\in E}-\psi^{y_0}(z)\de z}_{TV}\\
    &\esp\loq\dabs{\psf_{m_k}\pth{\widehat{Z}_k\in \de z\kt \widehat{\tau}>k\vg \widehat{Y}_k\in E}-\psf_{m}\pth{\widehat{Z}^{\eta}_k\in\de z\kt\widehat{\tau}^{\eta}>k\vg\widehat{Y}^{\eta}_k\in E}}_{TV}\\
    &\esp\esp+\dabs{\psf_{m}\pth{\widehat{Z}^{\eta}_k\in\de z\kt\widehat{\tau}^{\eta}>k\vg\widehat{Y}^{\eta}_k\in E}-\psi^{y_0}(z)\de z}_{TV}\\
    &\esp\loq K\pth{e^{\chi'\eta k}-1}+Ke^{-\chi k},
  \end{align*}
  for $K>0$ large enough.
  Hence,
  \begin{align*}
    &\dabs{\psf_{\varphi}\pth{Y_t\in E,Z_t\in \de z\kt \tau>t}-\psf_{\varphi}(Y_t\in E\kt \tau>t)\psi^{y_0}(z)\de z}_{TV}\\
    &\esp\loq K\sum_{k=0}^{\floor{t}}\psf_{\varphi}\pth{Y_t\in E,T(\eta;t)\in I_k\kt \tau>t}\pth{e^{-\chi k}+e^{\chi'\eta k}-1}\\
    &\esp\loq K\psf_{\varphi}\pth{Y_t\in E\kt \tau>t}\sum_{k=0}^{\floor{t}}\psf_{\varphi}\pth{T(\eta;t)\in I_k\kt \tau>t, Y_t\in E}\times e^{-\chi k}\\
    &\esp\esp+K\psf_{\varphi}\pth{Y_t\in E\kt \tau>t}\times\pth{e^{\chi'\eta t}-1}.
  \end{align*}
 This gives, approximating the sum from above by an expectation:
  \begin{align*}
    &\dabs{\psf_{\varphi}\pth{Y_t\in E,Z_t\in \de z\kt \tau>t}-\psf_{\varphi}(Y_t\in E\kt \tau>t)\psi^{y_0}(z)\de z}_{TV}\\    &\esp\loq K\psf_{\varphi}\pth{Y_t\in E\kt \tau>t}\esf_{\varphi}\cro{e^{-\chi (T(\eta;t)-1)}\kt \tau>t,Y_t\in E}\\
    &\esp\esp+K\psf_{\varphi}\pth{Y_t\in E\kt \tau>t}\times\pth{e^{\chi'\eta t}-1}.
  \end{align*}
  Dividing both sides by $\psf_{\varphi}\pth{Y_t\in E\kt \tau>t}$, we conclude:
  \begin{align*}
    &\dabs{\psf_{\varphi}\pth{Z_t\in \de z\kt \tau>t\vg Y_t\in E}-\psi^{y_0}(z)\de z}_{TV}\\
    &\esp\loq (Ke^{\chi})\esf_{\varphi}\cro{e^{-\chi T(\eta;t)}\kt \tau>t,Y_t\in E}+K\pth{e^{\chi'\eta t}-1}.
  \end{align*}
Finally, the statements of Lemmas~\ref{lem:cv_shat} and~\ref{lem:comparison_frozen} show that the constants $K$, $\chi$ and $\chi'$ can be chosen independently of $y_0$.
\end{proof}

\subsection{Conclusion of the proof}\label{ss:conclusion_big_thm}
In Subsections~\ref{ss:frozen} and~\ref{ss:comparison_frozen}, we worked with a fixed $\eps>0$ and omitted the superindex $\eps$. Here, we will make $\eps$ go to zero. Therefore, we need to go back to our former notations: $Y^{\eps}$, $Z^{\eps}$, $\tau^{\eps}$, $T(\eta,\eps;t)$, $\varphi_{\eps}$, etc.
First, we prove Proposition~\ref{ppn:estimation_diff_zero}. Then, we will be able to prove Theorem~\ref{thm:diffusion_zero}.

\begin{proof}[Proof of Proposition~\ref{ppn:estimation_diff_zero}]
    Let $\eta>0$. Recall that $V^{\eta}=\mathcal{B}(y_0,\eta)\cap\mathcal{Y}$. By Lemma~\ref{lem:estimation_diff_zero}, there exist $K>0$, $\chi>0$ and $\chi'>0$ such that for all $E\subset V^{\eta}$, for all $t>0$, $\eta>0$, $\eps>0$,
\begin{align}
  &\dabs{\psf_{\varphi_{\eps}}\pth{{Z}^{\eps}_t\in \de z\ |\ {\tau}^{\eps}>t,\, Y_t^{\eps}\in E}-\psi^{y_0}(z)\de z}_{TV}\nonumber\\
  &\esp\esp<K\esf_{\varphi_{\eps}}\cro{e^{-\chi T(\eta,\eps;t)}\ |\  \tau^{\eps}>t,\, Y^{\eps}_t\in E}+K\pth{e^{\eta\chi' t}-1}.\label{eq:csq_lemma}
\end{align}
For $\eps\in(0,1)$ sufficiently small, we let $t_{\eps}:=\frac{\ln(1/\eps)}{1+4(\max d-\min d)}$, $\eta_{\eps}:=\eps^{1/4}$, and $E^{\eps}:=B\pth{y_0,\eps^2}\cap\mathcal{Y}\subset V^{\eta_{\eps}}$.
  For~$\eps>0$ sufficiently small, we have $d(E^{\eps},\dr V^{\eta_{\eps}})\goq\eta_{\eps}/2$.
  Further, $\inf_{\eps\in(0,1)}\frac{\abs{E^{\eps}}}{\eps^2}>0$.
Our goal is to show that as $\eps\to 0$,
\begin{align}
  \esf_{\varphi_{\eps}}\cro{e^{-\chi T(\eta_{\eps},\eps;t_{\eps})}\ |\  \tau^{\eps}>t_{\eps},\, Y^{\eps}_{t_{\eps}}\in E^{\eps}}&\to 0,\label{eq:condition1}\\
  {e^{\eta_{\eps}\chi' t_{\eps}}-1}&\to 0.\label{eq:condition2}
\end{align}
With~\eqref{eq:csq_lemma}, we will then be able to conclude. 
First, we have $t_{\eps}\eta_{\eps}\to 0$ as $\eps\to 0$, so~\eqref{eq:condition2} holds. Second, let us show that~\eqref{eq:condition1} holds.
%Fix $t'>0$ and consider $\eps\in(0,1)$ small enough that $t_{\eps}>t'$.
We have
\begin{align}
  &\esf_{\varphi_{\eps}}\cro{e^{-\chi T(\eta_{\eps},\eps;t_{\eps})}\ |\  \tau^{\eps}>t_{\eps},\, Y^{\eps}_{t_{\eps}}\in E^{\eps}}\nonumber\\
  &\esp\loq e^{-\chi t_{\eps}}+\psf_{\varphi_{\eps}}\pth{T(\eta_{\eps},\eps;t_{\eps})< t_{\eps}\ |\  \tau^{\eps}>t_{\eps},\, Y^{\eps}_{t_{\eps}}\in E^{\eps}}.\label{eq:ppn_estimate_exp}
\end{align} %%%%%%%%%%%%
  We wish to estimate the very last term. Let $\widetilde{E}^{\eps}:=\acc{y\in\mathcal{Y}\ /\ d(y,E^{\eps})<\eta_{\eps}/4}$. Then $E^{\eps}\subset\widetilde{E}^{\eps}$ and the triangle inequality implies that for~$\eps>0$ sufficiently small, $d(\widetilde{E}^{\eps},\dr V^{\eta_{\eps}})\goq\eta_{\eps}/4$.
  We have, decomposing the above probability according to the value of $Y^{\eps}_0$,
\begin{align}
  &\psf_{\varphi_{\eps}}\pth{T(\eta_{\eps},\eps;t_{\eps})<t_{\eps}\ |\  \tau^{\eps}>t_{\eps},\, Y^{\eps}_{t_{\eps}}\in E^{\eps}}\nonumber\\
  &\esp\loq \psf_{\varphi_{\eps}}\pth{Y^{\eps}_0\notin\widetilde{E}^{\eps}\ |\  \tau^{\eps}>t_{\eps},\, Y^{\eps}_{t_{\eps}}\in E^{\eps}}+\psf_{\varphi_{\eps}}\pth{T(\eta_{\eps},\eps;t_{\eps})<t_{\eps},\,Y^{\eps}_0\in\widetilde{E}^{\eps}\ |\  \tau^{\eps}>t_{\eps},\, Y^{\eps}_{t_{\eps}}\in E^{\eps}}\nonumber\\
  &\esp\loq \psf_{\varphi_{\eps}}\pth{Y^{\eps}_0\notin\widetilde{E}^{\eps}\ |\  \tau^{\eps}>t_{\eps},\, Y^{\eps}_{t_{\eps}}\in E^{\eps}}+\psf_{\varphi_{\eps}}\pth{T(\eta_{\eps},\eps;t_{\eps})<t_{\eps}\ |\  \tau^{\eps}>t_{\eps},\, Y^{\eps}_{t_{\eps}}\in E^{\eps},\,Y^{\eps}_0\in\widetilde{E}^{\eps}}.\label{eq:estimate_main_ppn}
\end{align}

\paragraph{Step 1. Estimate the first term of~\eqref{eq:estimate_main_ppn}.}
First, let us show that the first term on the right-hand side of~\eqref{eq:estimate_main_ppn} converges to~$0$ as $\eps\to0$.
Since $d(\widetilde{E}^{\eps},\dr V^{\eta_{\eps}})\goq\eta_{\eps}/4$ for $\eps$ sufficiently small, we have:
\begin{align*}
  \psf_{\varphi_{\eps}}\pth{Y^{\eps}_0\notin\widetilde{E}^{\eps}\ |\  \tau^{\eps}>t_{\eps},\, Y^{\eps}_{t_{\eps}}\in E^{\eps}}
  &\loq \psf_{\varphi_{\eps}}\pth{\dabs{Y^{\eps}_0-Y^{\eps}_{t_{\eps}}}>\frac{\eta_{\eps}}{4}\ |\  \tau^{\eps}>t_{\eps},\, Y^{\eps}_{t_{\eps}}\in E^{\eps}}.
\end{align*}
%We have:
%\[e^{-(\max d)t_{\eps}}\loq\psf_{\varphi_{\eps}}\pth{\tau^{\eps}>t_{\eps}}\loq e^{-(\min d)t_{\eps}}.\]
Since the Poisson point process defining $\tau_{\eps}$ is independent of the carrying process $(\widetilde{Y}^{\eps}_t,\widetilde{Z}^{\eps}_t)_t$, there exist two exponential random variables $\underline{\tau}$ and $\overline{\tau}$, with rates $\max d$ and $\min d$ respectively, independent of $(\widetilde{Y}^{\eps}_t,\widetilde{Z}^{\eps}_t)$, and such that
\begin{equation}\label{eq:coupling_tau}
  \underline{\tau}\loq\tau_{\eps}\loq\overline{\tau}.
\end{equation}
Using this comparison,
we have:
\begin{align*}
  &\psf_{\varphi_{\eps}}\pth{\dabs{Y^{\eps}_0-Y^{\eps}_{t_{\eps}}}>\frac{\eta_{\eps}}{4}\ |\  \tau^{\eps}>t_{\eps},\, Y^{\eps}_{t_{\eps}}\in E^{\eps}}\\
  &\qquad\loq\frac{1}{\psf_{\varphi_{\eps}}\pth{\underline{\tau}>t_{\eps},\, \widetilde{Y}^{\eps}_{t_{\eps}}\in E^{\eps}}}\psf_{\varphi_{\eps}}\pth{\dabs{\widetilde{Y}^{\eps}_0-\widetilde{Y}^{\eps}_{t_{\eps}}}>\frac{\eta_{\eps}}{4},\,\overline{\tau}>t_{\eps},\, \widetilde{Y}^{\eps}_{t_{\eps}}\in E^{\eps}}\\
  &\qquad\loq\psf_{\varphi_{\eps}}\pth{\dabs{\widetilde{Y}^{\eps}_0-\widetilde{Y}^{\eps}_{t_{\eps}}}>\frac{\eta_{\eps}}{4}\ |\  \widetilde{Y}^{\eps}_{t_{\eps}}\in E^{\eps}}e^{(\max d-\min d)t_{\eps}}\\
  &\qquad\loq\frac{\psf_{\varphi_{\eps}}\pth{\dabs{\widetilde{Y}^{\eps}_0-\widetilde{Y}^{\eps}_{t_{\eps}}}>\frac{\eta_{\eps}}{4}}}{\psf_{\varphi_{\eps}}\pth{\widetilde{Y}^{\eps}_{t_{\eps}}\in E^{\eps}}}e^{(\max d-\min d)t_{\eps}}.
%  &\qquad\loq \frac{16}{\eta_{\eps}^2}\,\mathds{E}_{\varphi_{\eps}}\cro{\dabs{Y^{\eps}_0-Y^{\eps}_{t_{\eps}}}^2\ |\  Y^{\eps}_{t_{\eps}}\in E^{\eps}}e^{(\max d-\min d)t_{\eps}},
\end{align*}
%where we used the Markov inequality for the last line.
%Hence, there exists a constant $C>0$ depending only on the diffusion coefficient~$A$ and the drift term~$B$, but not on~$\eps$ or~$y_0$, such that %(using the Markov inequality for the first line)
%\begin{align*}
%  \psf_{\varphi_{\eps}}\pth{\dabs{Y^{\eps}_0-Y^{\eps}_{t_{\eps}}}>\frac{\eta_{\eps}}{4}\ |\  \tau^{\eps}>t_{\eps},\, Y^{\eps}_{t_{\eps}}\in E^{\eps}}
%  %&\loq \frac{16}{\eta_{\eps}^2}\,\mathds{E}_{\varphi_{\eps}}\cro{\dabs{Y^{\eps}_0-Y^{\eps}_{t_{\eps}}}^2\ |\  Y^{\eps}_{t_{\eps}}\in E^{\eps}}e^{(\max d-\min d)t_{\eps}}\\
%  &\loq \frac{C\eps^2 t_{\eps}}{\eta_{\eps}^2}e^{(\max d-\min d)t_{\eps}}\\
%  &= 4C\eps t_{\eps}e^{(\max d-\min d)t_{\eps}}.
%\end{align*}
On the one hand, there exist constants $C_1,C_2>0$, depending only on the diffusion coefficient~$A$, the drift term~$B$ and the domain~$\mathcal{Y}$ (but not on~$\eps$ or~$y_0$), such that
\[\psf_{\varphi_{\eps}}\pth{\dabs{\widetilde{Y}^{\eps}_0-\widetilde{Y}^{\eps}_{t_{\eps}}}>\frac{\eta_{\eps}}{4}}\loq C_1e^{-C_2\eta_{\eps}^2/\eps^2}=C_1e^{-C_2/\eps^{3/2}}.\]
On the other hand, let us estimate $\psf_{\varphi_{\eps}}\pth{\widetilde{Y}^{\eps}_{t_{\eps}}\in E^{\eps}}$.
The Neumann boundary condition on~$\dr\mathcal{Y}$ or~$\dr\mathcal{Z}$ (if nonempty) allows us to apply the Harnack inequality up to the boundary to the function $(y,z)\mapsto\varphi_{\eps}(y/\eps,z)$,
which solves an equation that is uniformly elliptic in~$\eps$ and that has coefficients uniformly bounded in~$\eps$ (for more details, see the beginning of the proof of Theorem~\ref{thm:diffusion_zero}).
We obtain that there exists a constant $C_3\in(0,1)$, depending only on~$A$,~$B$,~$r$, but not on~$\eps$ or~$y_0$,
such that for all $y\in\mathcal{Y}$ (denoting by $\mathcal{Y}[y,\eps]$ the connected component of~$\mathcal{Y}\cap B(y,\eps)$ containing~$y$),
\begin{equation}\label{eq:loc_harnack_proof_prop}
  \inf_{y'\in\mathcal{Y}[y,\eps]}\ \varphi_{\eps}(y',\mathcal{Z})>C_3\varphi_{\eps}(y,\mathcal{Z}).
\end{equation}
Let us extend this inequality to any pair of points $y,y'\in\mathcal{Y}$.
Since~$\mathcal{Y}$ is connected and smooth,~$\overline{\mathcal{Y}}$ is path-connected, with paths that can be chosen to be~$\mcc^1$. Hence we can define the geodesic distance $d_{\overline{\mcy}}$ on~$\overline{\mcy}$. % between $y\in\overline{\mathcal{Y}}$ and $y'\in\overline{\mathcal{Y}}$ as
%\[d_{\overline{\mathcal{Y}}}(y,y'):=\inf\acc{\pth{\int_0^1\abs{\nabla\gamma}^2}^{1/2}\ /\ \gamma\in\mathcal{C}^1([0,1],\overline{\mathcal{Y}}),\ \gamma(0)=y,\ \gamma(1)=y'}.\]
By continuity of~$d_{\overline{\mathcal{Y}}}$ and compactness of~$\overline{\mathcal{Y}}$, the radius $R:=\sup_{y,y'\in{\overline{\mathcal{Y}}}}\acc{d_{\mathcal{Y}}(y,y')}$ is finite.
Hence, for any $y,y'\in\mathcal{Y}$, there exists a sequence $y=y_1,y_2,\hdots,y_k=y'$, with $k\loq{R/\eps}+1$, such that $\mathcal{Y}[y_i,\eps]\cap\mathcal{Y}[y_{i+1},\eps]$ is nonempty for all $i=1,\hdots,k-1$.
Therefore, we may apply~\eqref{eq:loc_harnack_proof_prop} (at most) $\ceil{R/\eps}+1$ times to obtain: for all $\eps>0$,
\[\inf_{y\in\mathcal{Y}}\varphi_{\eps}(y,\mathcal{Z})>C_3^{\ceil{R/\eps}+1}\sup_{y\in\mathcal{Y}}\varphi_{\eps}(y,\mathcal{Z}).\]
Recalling the normalisation~$\varphi_{\eps}(\mathcal{Y},\mathcal{Z})=1$, we have $\sup_{y\in\mathcal{Y}}\varphi_{\eps}(y,\mathcal{Z})\goq 1/\abs{\mathcal{Y}}$, so
\[\psf_{\varphi_{\eps}}\pth{\widetilde{Y}^{\eps}_{0}\in E^{\eps}}\goq \abs{E^{\eps}}\inf_{y\in\mathcal{Y}}\varphi_{\eps}(y,\mathcal{Z})\goq \frac{\abs{E^{\eps}}C_3^{\ceil{R/\eps}+1}}{\abs{\mathcal{Y}}},\]
which implies:
\[\psf_{\varphi_{\eps}}\pth{\widetilde{Y}^{\eps}_{t_{\eps}}\in E^{\eps}}\goq
\pth{\inf_{y\in E^{\eps}}\psf_{y}\pth{\widetilde{Y}^{\eps}_{1}\in E^{\eps}}}^{\ceil{t_\eps}}\psf_{\varphi_{\eps}}\pth{\widetilde{Y}^{\eps}_{0}\in E^{\eps}}
\goq
\frac{\abs{E^{\eps}}C_3^{\ceil{R/\eps}+1}(\eps^2 C_4)^{\ceil{t_{\eps}}}}{\abs{\mathcal{Y}}},\]
with
\[C_4:=\inf_{\eps\in(0,1),\, y_0'\in\mathcal{Y}}\ \cro{\frac{1}{\eps^2}\inf_{y\in B(y_0',\eps^2)\cap\mathcal{Y}}\psf_{y}\pth{\widetilde{Y}^{\eps}_{1}\in B(y_0',\eps^2)\cap\mathcal{Y}}}>0,\]
which depends only on the diffusion coefficient~$A$ and the drift coefficient~$B$, but not on~$\eps$ or~$y_0$ (note that as $\eps\to0$, the term between the brackets converges to~$+\infty$ uniformly in~$y_0'$, so~$C_4$ is indeed positive).
We obtain:
\begin{align*}
  \psf_{\varphi_{\eps}}\pth{\widetilde{Y}^{\eps}_0\notin\widetilde{E}^{\eps}\ |\  \tau^{\eps}>t_{\eps},\, \widetilde{Y}^{\eps}_{t_{\eps}}\in E^{\eps}}
  &\loq\frac{\abs{\mathcal{Y}}C_1e^{-C_2/\eps^{3/2}}}{\abs{E^{\eps}}C_3^{\ceil{R/\eps}+1}(\eps^2C_4)^{\ceil{t_{\eps}}}}e^{(\max d-\min d)t_{\eps}}.
\end{align*}
This gives, recalling that $\inf_{\eps\in(0,1)}\abs{E^{\eps}}/\eps^2>0$,
\begin{equation}\label{eq:estimate_main_ppn_left}
  \psf_{\varphi_{\eps}}\pth{\widetilde{Y}^{\eps}_0\notin\widetilde{E}^{\eps}\ |\  \tau^{\eps}>t_{\eps},\, \widetilde{Y}^{\eps}_{t_{\eps}}\in E^{\eps}}\ \mathop{\longrightarrow}_{\eps\to0}\ 0.
\end{equation}

\paragraph{Step 2. Estimate the second term of~\eqref{eq:estimate_main_ppn} and conclude.}
Second, let us show that the second term on the right-hand side of~\eqref{eq:estimate_main_ppn} converges to~$0$ as $\eps\to0$.
Define the stopping time
\[\tau_V:=\inf\acc{u>0\ / \ \widetilde{Y}^{\eps}_u\in\dr V^{\eta_{\eps}}}.\]
On the event $\acc{T(\eta_{\eps},\eps;t_{\eps})< t_{\eps}}\cap\acc{\widetilde{Y}^{\eps}_{0}\in \widetilde{E}^{\eps},\,\widetilde{Y}^{\eps}_{t_{\eps}}\in {E}^{\eps}}$, we have $0<\tau_V<t_{\eps}$.
Therefore, 
\begin{align*}
  &\psf_{\varphi_{\eps}}\pth{T(\eta_{\eps},\eps;t_{\eps})<t_{\eps}\ |\  \tau^{\eps}>t_{\eps},\, \widetilde{Y}^{\eps}_{t_{\eps}}\in E^{\eps},\,\widetilde{Y}^{\eps}_0\in\widetilde{E}^{\eps}}\\
  &\qquad\loq\psf_{\varphi_{\eps}}\pth{\tau_V\in(0,t_{\eps})\ |\  \tau^{\eps}>t_{\eps},\, \widetilde{Y}^{\eps}_{t_{\eps}}\in E^{\eps},\,\widetilde{Y}^{\eps}_0\in\widetilde{E}^{\eps}}\\
  &\qquad\loq \psf_{\varphi_{\eps}}\pth{\tau_V\in(0,t_{\eps})\ |\  \widetilde{Y}^{\eps}_{t_{\eps}}\in E^{\eps},\,\widetilde{Y}^{\eps}_0\in\widetilde{E}^{\eps}}\,e^{(\max d-\min d)t_{\eps}}.
\end{align*}
The last line is obtained by using, as above, the comparison~\eqref{eq:coupling_tau}.
Further, since $E^{\eps}\subset\widetilde{E}^{\eps}\subset V^{\eta_{\eps}}$, we have:
%\[\psf_{\varphi_{\eps}}\pth{\widetilde{Y}^{\eps}_{t_{\eps}}\in E^{\eps} \ |\ \tau_V\in(0,t_{\eps}),\, \widetilde{Y}^{\eps}_0\in\widetilde{E}^{\eps}}
%\loq \psf_{\varphi_{\eps}}\pth{\widetilde{Y}^{\eps}_{t_{\eps}}\in E^{\eps} \ |\ \widetilde{Y}^{\eps}_0\in\widetilde{E}^{\eps}},\]
%which implies:
\[\psf_{\varphi_{\eps}}\pth{\tau_V\in(0,t_{\eps})\ |\  \widetilde{Y}^{\eps}_{t_{\eps}}\in E^{\eps},\,\widetilde{Y}^{\eps}_0\in\widetilde{E}^{\eps}}
\loq \psf_{\varphi_{\eps}}\pth{\tau_V\in(0,t_{\eps})\ |\  \widetilde{Y}^{\eps}_0\in\widetilde{E}^{\eps}}.
\]
%Further, by the strong Markov property, conditionnally on~$\widetilde{Y}^{\eps}_{\tau_V}$, $(\widetilde{Y}^{\eps}_{t})_{t\goq \tau_V}$ is independent of~$\widetilde{Y}^{\eps}_0$.
%Therefore,
%\begin{align*}
%  \psf_{\varphi_{\eps}}\pth{\tau_V\in(0,t_{\eps})\ |\  \widetilde{Y}^{\eps}_{t_{\eps}}\in E^{\eps},\,\widetilde{Y}^{\eps}_0\in\widetilde{E}^{\eps}}
%  &=\int_{\dr V^{\eta}}\psf_{\varphi_{\eps}}\pth{\tau_V\in(0,t_{\eps}),\,\widetilde{Y}^{\eps}_{\tau_V}\in \de y\ |\  \widetilde{Y}^{\eps}_{t_{\eps}}\in E^{\eps},\,\widetilde{Y}^{\eps}_0\in\widetilde{E}^{\eps}}\\
%  &=\int_{\dr V^{\eta}}\psf_{\varphi_{\eps}}\pth{\tau_V\in(0,t_{\eps}),\,\widetilde{Y}^{\eps}_{\tau_V}\in \de y\ |\ \widetilde{Y}^{\eps}_0\in\widetilde{E}^{\eps}}\\
%  &=
%   & \loq \psf_{\varphi_{\eps}}\pth{\tau_V\in(0,t_{\eps})\ |\  \widetilde{Y}^{\eps}_0\in\widetilde{E}^{\eps}}.
%\end{align*}
Last, recall that $d(\widetilde{E}^{\eps},\dr V^{\eta})\goq \eta_{\eps}/4$ for $\eps>0$ sufficiently small.
Thus, if $\widetilde{Y}^{\eps}_0\in\widetilde{E}^{\eps}$, then $\|\widetilde{Y}^{\eps}_{\tau_V}-\widetilde{Y}^{\eps}_0\|\goq\frac{\eta_{\eps}}{4}$.
Therefore, using the Markov inequality,
\begin{align*}
  &\psf_{\varphi_{\eps}}\pth{\tau_V\in(0,t_{\eps})\ |\  \widetilde{Y}^{\eps}_0\in\widetilde{E}^{\eps}}\\
  &\qquad\loq\psf_{\varphi_{\eps}}\pth{\sup_{u\in[0,t_{\eps}]}\|\widetilde{Y}^{\eps}_{u}-\widetilde{Y}^{\eps}_0\|\goq\frac{\eta_{\eps}}{4} |\  \widetilde{Y}^{\eps}_0\in\widetilde{E}^{\eps}}
\loq  \frac{C''t_{\eps}\eps^2}{\eta_{\eps}^2}=C''t_{\eps}\eps^{3/2}, %\label{eq:estimate_up_ppn}
\end{align*}
for some constant $C''>0$ depending on the diffusion coefficient~$A$ and the drift coefficient~$B$, but not depending on~$\eps$ or~$y_0$.
Combining this with the above, we obtain:
\begin{align*}
  \psf_{\varphi_{\eps}}\pth{T(\eta_{\eps},\eps;t_{\eps})< t_{\eps}\ |\  \tau^{\eps}>t_{\eps},\, Y^{\eps}_{t_{\eps}}\in E^{\eps},\,Y^{\eps}_0\in\widetilde{E}^{\eps}}
  &\loq C''t_{\eps}\eps^{3/2}\,e^{(\max d-\min d)t_{\eps}},
\end{align*}
which converges to~$0$ as $\eps\to0$.
Therefore~\eqref{eq:estimate_main_ppn} and~\eqref{eq:estimate_main_ppn_left} imply that as $\eps\to0$,
\[  \psf_{\varphi_{\eps}}\pth{T(\eta_{\eps},\eps;t_{\eps})<t_{\eps}\ |\  \tau^{\eps}>t_{\eps},\, Y^{\eps}_{t_{\eps}}\in E^{\eps}}\to 0.\]
Thus, taking $\eps\to0$ in Equation~\eqref{eq:ppn_estimate_exp} implies that~\eqref{eq:condition1} holds.
The convergence stated in the proposition follows from ~\eqref{eq:csq_lemma} combined with~\eqref{eq:condition1}--\eqref{eq:condition2}.
Since the bound of Lemma~\ref{lem:estimation_diff_zero} is uniform in $y_0\in\overline{\mathcal{Y}}$, we conclude that the convergence is also uniform in $y_0\in\overline{\mathcal{Y}}$.
\end{proof}

We are ready to prove Theorem~\ref{thm:diffusion_zero}, using Proposition~\ref{ppn:estimation_diff_zero}.

\begin{proof}[Proof of Theorem~\ref{thm:diffusion_zero}]
  For $\eps\in\ioo{0,1}$, we set $\eta_{\eps}:=\eps^{1/4}$ and we let $E^{\eps}:=B\pth{y_0,\eps^2}\cap \mcy\subset V^{\eta_{\eps}}$.
  
  Let us change coordinates. For $\eps\in\ioo{0,1}$, we define \[\mcy_{\eps}:=\acc{y\in\Er^n\tq y_{0}+\eps y\in\mcy}.\]
  For all $y\in{\mcy_{\eps}}$, we set
  \[\widetilde{\varphi}_{\eps}(y,z):=\frac{\varphi_{\eps}\pth{y_0+\eps y
,z}}{\varphi_{\eps}\pth{y_0,\mcz}}.\]
We define $\widetilde{A}$, $\widetilde{B}$, $\widetilde{\mcl}_z$ and $\widetilde{c}$ with the same change of variables (we omit the dependence in $\eps$ for the sake of clarity).
We let $\widetilde{E}^{\eps}:=B\pth{0,\eps}\cap\mcy_{\eps}$ be the projection of $E^{\eps}$ in this new system of coordinates. 
  First, $\widetilde{\varphi}_{\eps}$ satisfies
  \begin{equation*}
  \nabla_y\cdot(\widetilde{A}\nabla_y\widetilde{\varphi}_\eps)+\nabla_y\cdot\pth{\widetilde{\varphi}_\eps \widetilde{B}}+(\widetilde{\mcl}_z+\widetilde{c})\widetilde{\varphi}_\eps=k_\eps\widetilde{\varphi}_{\eps},\esp y\in\mcy_{\eps},\esp z\in\mcz.  
  \end{equation*}
  Thanks to this transformation, the diffusion coefficients become uniformly elliptic on the whole domain. Let $U\subset\subset \mcy_1$, so that $U\subset\subset\mcy_{\eps}$ for all $\eps\in\ioo{0,1}$. Let $\mcz'\subset\subset\mcz$. For all $\eps>0$, we have $\widetilde{\varphi}_{\eps}(0,\mcz)= 1$.
  Therefore, by the Harnack inequality, the family $(\widetilde{\varphi}_{\eps})_{\eps\in\ioo{0,1}}$ is bounded in $L^{\infty}(U\times\mcz')$.
  By the interior Schauder inequality (\cite[Theorem 6.2]{GilTru01}), therefore, the family $(\widetilde{\varphi}_{\eps})_{\eps\in\ioo{0,1}}$ is bounded in $\mcc^{2,\alpha}(U\times\mcz')$. Hence, taking $\beta\in\ioo{0,\alpha}$, up to extraction of a subsequence, there exists a limit $\widetilde{\varphi}_0\in \mcc^{2,\beta}(U\times\mcz')$:
\[\widetilde{\varphi}_{\eps}\to \widetilde{\varphi}_0\comment{in $\mcc^{2,\beta}(U\times\mcz')$}.\]
Since $\abs{\widetilde{E}^{\eps}}\to 0$, this implies in particular:
\begin{equation}\label{eq:cv_varphitilde_y0}
  \frac{1}{\abs{\widetilde{E}^{\eps}}}\widetilde{\varphi}_{\eps}(\widetilde{E}^{\eps},\cdot)=\frac{1}{\abs{\widetilde{E}^{\eps}}}\int_{\widetilde{E}^{\eps}}\widetilde{\varphi}_{\eps}(y,\cdot)\de y\to \widetilde{\varphi}_0(0,\cdot)\comment{in $\mcc^{0}(\mcz')$}.
\end{equation}
%Recall that $\widetilde{\varphi}_{\eps}(\widetilde{E}^{\eps},\cdot)$ is an abuse of notation to denote \[\widetilde{\varphi}_{\eps}(\widetilde{E}^{\eps},\cdot)=\int_{\widetilde{E}^{\eps}}\widetilde{\varphi}_{\eps}(y,\cdot)\de y.\]
Finally, we note that
\[\frac{1}{\abs{\widetilde{E}^{\eps}}}\widetilde{\varphi}_{\eps}(\widetilde{E}^{\eps},\cdot)
=\frac{\widetilde{\varphi}_{\eps}(\widetilde{E}^{\eps},\mcz)}{{\abs{\widetilde{E}^{\eps}}}}
\times\frac{\varphi_{\eps}(E^{\eps},\cdot)}{\varphi_{\eps}(E^{\eps},\mcz)}.\]
On the one hand,
\[\frac{\widetilde{\varphi}_{\eps}(\widetilde{E}^{\eps},\mcz)}{{\abs{\widetilde{E}^{\eps}}}}\to{\widetilde{\varphi}_{0}(0,\mcz)}=1.\]
On the other hand, using the notations of Proposition~\ref{ppn:estimation_diff_zero},
  \begin{equation*}
    \frac{\varphi_{\eps}(E^{\eps},\cdot)}{\varphi_{\eps}(E^{\eps},\mcz)}
    =\frac{\psf_{\varphi_\eps}\pth{Y^{\eps}_{t_{\eps}}\in E^{\eps},\, {Z}^{\eps}_{t_{\eps}}\in  \cdot\kt \tau^{\eps}>t_{\eps}}}
    {\psf_{\varphi_\eps}\pth{Y^{\eps}_{t_{\eps}}\in E^{\eps}\kt \tau^{\eps}>t_{\eps}}}
    =\psf_{\varphi_\eps}\pth{{Z}^{\eps}_{t_{\eps}}\in  \cdot\kt \tau^{\eps}>t_{\eps}\vg Y_{t_{\eps}}^{\eps}\in E^{\eps}}.
  \end{equation*}
Therefore, by Proposition~\ref{ppn:estimation_diff_zero},
  \begin{equation*}
  \frac{\varphi_{\eps}(E^{\eps},\cdot)}{\varphi_{\eps}(E^{\eps},\mcz)}\,\de z\to \psi^{y_0}(z)\de z
  \end{equation*}
  in total variation on $\mcz$.
Therefore, we have in total variation on $\mcz$ (thus also on~$\mcz'$):
\begin{equation}\label{eq:cv_varphitilde_tv}
  \frac{1}{\abs{\widetilde{E}^{\eps}}}\widetilde{\varphi}_{\eps}(\widetilde{E}^{\eps},z)\de z\to\psi^{y_0}(z)\de z\comment{as $\eps\to 0$}.
\end{equation}
Note that convergence in $\mcc^{2,\beta}(\mcz')$ of the density implies that the corresponding measure converges %implies convergence
in total variation on $\mcz'$. Combining~\eqref{eq:cv_varphitilde_tv} with~\eqref{eq:cv_varphitilde_y0} yields: for all $z\in\mcz'$,
\[\widetilde{\varphi}_{0}(0,z)=\psi^{y_0}(z).\]
Hence, as $\eps\to 0$,
\[\frac{{\varphi}_{\eps}(y_0,\cdot)}{{\varphi}_{\eps}(y_0,\mcz)}\to\psi^{y_0}\comment{in $\mcc^{2,\beta}(\mcz')$}.\]
This concludes the first point of the theorem. Finally, the convergence in total variation on $\mcz$ is uniform in $y_0\in\overline{\mcy}$, see Proposition~\ref{ppn:estimation_diff_zero}. Thus the second point holds.
\end{proof}

\section{Behaviour of the principal eigenvalue (Theorems~\ref{thm:bzero} and~\ref{thm:general})}\label{s:proof_cor}

In the previous section, we worked on the principal eigenfunction and proved Theorem~\ref{thm:diffusion_zero} using a probabilistic interpretation.
  We now use Theorem~\ref{thm:diffusion_zero} to prove results about the principal eigenvalue, in some particular cases.
We drop the probabilistic interpretation; the techniques in this part will be based on Hamilton-Jacobi equations.

\subsection{Preliminary results}

We record several lemmas which will be useful for the proof of Theorems~\ref{thm:bzero} and~\ref{thm:general}.
The first lemma is basic.
  The following two lemmas add technical precision on the convergence of $\varphi_{\eps}$, and follow from Theorem~\ref{thm:diffusion_zero}.
  The last lemma, which is at the heart of the proof of Theorems~\ref{thm:bzero} and~\ref{thm:general}, makes a connection between the limiting principal eigenvalue and a Hamilton-Jacobi equation.
  Studying this Hamilton-Jacobi equation will allow us to prove Theorems~\ref{thm:bzero} and~\ref{thm:general}.

\begin{lem}\label{lem:estimate_k}
  For all $\eps>0$,
  \[\abs{k_{\eps}}\loq\dabs{c}_{\infty}.\]
\end{lem}

\begin{proof}
  For $g\in \mcc^{0,\alpha}(\mcy\times\mcz)$, let $\kappa_{\eps}(g)$ be the principal eigenvalue associated with the operator
    \[\mcl_g:\phi\mapsto\eps^2\nabla_y\cdot(A\nabla_y\phi)+\eps\nabla_y\cdot(\phi B)+\pth{\mcl_z+g}\phi,\]
    %\[\mcl_g:\phi\mapsto\eps^2A\Delta_y\phi+\eps B\cdot\nabla_y\phi+\pth{\mcl_z+g}\phi,\]
    with Neumann boundary conditions on $\dr(\mcy\times\mcz)$. By definition, we have $k_\eps=\kappa_{\eps}(c)$.
  Moreover, if $g_1\loq g_2$ then $\kappa_{\eps}(g_1)\loq \kappa_{\eps}(g_2)$. Therefore,
  \[\kappa_{\eps}(-\dabs{c}_{\infty})\loq k_{\eps}\loq \kappa_{\eps}(\dabs{c}_{\infty}).\]
  Finally, when $g$ is constant, $\varphi\equiv 1$ is a principal eigenfunction of $\mcl_g$. Thus:
  \begin{align*}
    \kappa_{\eps}(-\dabs{c}_{\infty})=-\dabs{c}_{\infty},\esp \kappa_{\eps}(\dabs{c}_{\infty})=\dabs{c}_{\infty}.
  \end{align*}
  The conclusion follows.
\end{proof}

\begin{lem}\label{lem:bounded_family_yeps}
  Let $(y_\eps)_{\eps\in\ioo{0,1}}$ be a family of elements of $\overline{\mcy}$. For all $\mcz'\subset\subset\mcz$, the family
  \[\pth{\frac{\varphi_{\eps}(y_{\eps},\cdot)}{\varphi_{\eps}(y_{\eps},\mcz)}}_{\eps\in\ioo{0,1}}\]
  is bounded in $\mcc^{2,\alpha}(\mcz')$.
\end{lem}

\begin{proof}
  Let us change coordinates. For $\eps\in\ioo{0,1}$, we define $\mcy_{\eps}:=\acc{y\in\Er^n\tq y_{\eps}+\eps y\in\mcy}$.
  For all $y\in\overline{\mcy_{\eps}}$, we set
  \[\widetilde{\varphi}_{\eps}(y,z):=\frac{\varphi_{\eps}\pth{y_{\eps}+\eps y,z}}{\varphi_{\eps}\pth{y_{\eps},\mcz}},\]
  and we define $\widetilde{A}$, $\widetilde{B}$, $\widetilde{\mcl}_z$ and $\widetilde{c}$ with the same change of coordinates (we omit the dependence in $\eps$ for the sake of clarity).
  We have
  \begin{equation}\label{eq:varphi_tilde}
  \nabla_y\cdot(\widetilde{A}\nabla_y\widetilde{\varphi}_\eps)+\nabla_y\cdot\pth{\widetilde{\varphi}_\eps \widetilde{B}}+(\widetilde{\mcl}_z+\widetilde{c})\widetilde{\varphi}_\eps=k_\eps\widetilde{\varphi}_{\eps},\esp y\in\mcy_{\eps}\esp z\in\mcz.
  \end{equation}
  Thanks to this transformation, the diffusion coefficients become uniformly elliptic on the whole domain.
  
  Since $\dabs{\varphi_{\eps}(y_{\eps},\cdot)}_{\mcc^{2,\alpha}(\mcz')}=\dabs{\widetilde{\varphi}_{\eps}(0,\cdot)}_{\mcc^{2,\alpha}(\mcz')}$, it is enough for our scope to show that the family 
$(\widetilde{\varphi}_{\eps}(0,\cdot))_{\eps\in\ioo{0,1}}$ is bounded in~$\mcc^{2,\alpha}(\mcz')$.
A slight difficulty arises when $y_{\eps}$ is near the boundary of~$\mcy_{\eps}$.
We first deal with the particular case where 
\begin{equation}\label{eq:condition_far_bdy}
  \exists\delta>0,\ \forall\eps>0,\qquad d(y_{\eps},\dr\mcy_{\eps})>\delta.
\end{equation}

\paragraph{Case 1. Condition~\eqref{eq:condition_far_bdy} holds.}
Let $U:=B(0,\delta/2)$.
Then for all $\eps\in\ioo{0,1}$, we have $d(U,\dr\mcy_{\eps})>\delta/2$, so we can apply an interior Harnack inequality: the family $(\widetilde{\varphi}_{\eps})_{\eps\in\ioo{0,1}}$ is bounded in $L^{\infty}(U\times\mcz')$.
To conclude, we apply the interior Schauder inequalities (\cite[Theorem 6.2]{GilTru01}): the family $(\widetilde{\varphi}_{\eps})_{\eps\in\ioo{0,1}}$ is bounded in $\mcc^{2,\alpha}(U\times\mcz')$. Thus the family 
$(\widetilde{\varphi}_{\eps}(0,\cdot))_{\eps\in\ioo{0,1}}$ is bounded in~$\mcc^{2,\alpha}(\mcz')$.

\paragraph{Case 2. Condition~\eqref{eq:condition_far_bdy} does not hold.} Let $U_{\eps}:=B(0,1)\cap\mcy_{\eps}$ and let
\[I:=\acc{\eps\in\ioo{0,1}\tq d(y_{\eps},\mcy_{\eps})\loq 1}.\]
By Case~1, there remains to prove that $(\dabs{\widetilde{\varphi}_{\eps}}_{\mcc^{2,\alpha}(U_{\eps}\times\mcz')})_{\eps\in I}$ is bounded.
We use the fact that $\widetilde{\varphi}_{\eps}$ satisfies~\eqref{eq:varphi_tilde} and the boundary condition $\nu\cdot\mca\nabla\widetilde{\varphi}_{\eps}=0$. There exists a constant $C>0$ such that for all $\eps\in\ioo{0,1}$, the $\mcc^{2,\alpha}$-norm of $\dr\mcy_{\eps}$ is at most $C$, \ie for all $\eps>0$ and for all $y\in\dr\mcy_{\eps}$, there exists a diffeomorphism $\Phi$ which straightens $\dr\mcy_{\eps}$ near $y$ and which satisfies $\dabs{\Phi}_{\mcc^{2,\alpha}}\loq C$.

Thus we can apply the elliptic Harnack inequality up to the boundary $\dr\mcy_{\eps}$, so there exists a constant~$C'$ such that for all~$\eps>0$, we have $\dabs{\widetilde{\varphi}_{\eps}}_{L^{\infty}(U_{\eps}\times\mcz)}\loq C'$.
Then, by boundary Schauder estimates (\cite{GilTru01}, see the proof of Theorem 6.30 and the comments after the proof of Lemma 6.5), the family $(\dabs{\widetilde{\varphi}_{\eps}}_{\mcc^{2,\alpha}(U_{\eps}\times\mcz')})_{\eps\in I}$ is bounded.
Again, we conclude that the family 
$(\widetilde{\varphi}_{\eps}(0,\cdot))_{\eps\in\ioo{0,1}}$ is bounded in $\mcc^{2,\alpha}(\mcz')$.
\end{proof}

\begin{lem}\label{lem:cv_to_ky}
  Let $\eps_n>0$, $\eps_n\to 0$ and let $(y_n)_n$ be a sequence of elements $y_n\in\overline{\mcy}$ converging to some $y_{\infty}\in\overline{\mcy}$. Let $z_0\in\mcz$. Then, up to extraction of a subsequence,
  \[\frac{(\mcl_z+c)\varphi_{\eps_n}(y_n,z_0)}{\varphi_{\eps_n}(y_n,z_0)}\to k^{y_{\infty}}.\]
\end{lem}
\begin{proof}
  Take $\mcz'\subset\subset\mcz$ and $\beta\in\ioo{0,\alpha}$. We use the shortcut $\varphi_n:=\varphi_{\eps_n}$. Lemma~\ref{lem:bounded_family_yeps} implies that $\pth{\frac{\varphi_n(y_n,\cdot)}{\varphi_n(y_n,\mcz)}}_{n\goq 0}$ is bounded in $\mcc^{2,\alpha}(\mcz')$. Hence up to extraction we must have in $\mcc^{2,\beta}(\mcz')$:
  \begin{equation}\label{eq:cv_to_psi}
    \frac{\varphi_n(y_n,\cdot)}{\varphi_n(y_n,\mcz)}\to \psi
  \end{equation}
  for some $\psi\in\mcc^{2,\beta}(\mcz')$.
  Moreover, by Theorem~\ref{thm:diffusion_zero}, we have
    \[\frac{\varphi_n(y_n,z)}{\varphi_n(y_n,\mcz)}\de z\to \psi^{y_{\infty}}(z)\de z\]
  in total variation of measures (we crucially use here that the convergence holds uniformly in~$y$). This implies that~\eqref{eq:cv_to_psi}
  %~\eqref{eq:cv_to_psi} holds in total variation norm on $\mcz'$ with $\psi=\psi^{y_{\infty}}$.
  %Note that convergence in $\mcc^{2,\beta}(\mcz')$ implies convergence in total variation on~$\mcz'$. Thus
  holds in $\mcc^{2,\beta}(\mcz')$ with $\psi=\psi^{y_{\infty}}$.
 To conclude, we write:
\[\frac{(\mcl_z+c)\varphi_n(y_n,z_0)}{\varphi_n(y_n,z_0)}=\frac{(\mcl_z+c)\varphi_n(y_n,z_0)}{\varphi_n(y_n,\mcz)}\times\frac{\varphi_n(y_n,\mcz)}{\varphi_n(y_n,z_0)}.\]
Thus
\[\lim_{n\to+\infty}\frac{(\mcl_z+c)\varphi_n(y_n,z_0)}{\varphi_n(y_n,z_0)}=\frac{(\mcl_z+c)\psi^{y_{\infty}}(z_0)}{\psi^{y_{\infty}}(z_0)}=k^{y_{\infty}}.\]
\end{proof}

By Lemma~\ref{lem:estimate_k}, the family $(k_{\eps})_{\eps>0}$ is bounded. Thus, there exists a sequence $(\eps_n)_n$ such that $\eps_n\to 0$ and $(k_{\eps_n})_{n\goq 0}$ converges to a limit $k_0$. The following lemma, inspired by Proposition 3.1 in~\cite{HNR11}, gives a property of $k_0$.
This property is crucial for the proofs of Theorems~\ref{thm:bzero} and~\ref{thm:general}.

\begin{lem}\label{lem:visco_sol}
  Let $k_0\in\Er$ and $\mcy'\subset\subset\mcy$. If there exists a sequence $(\eps_n)_{n\goq 0}$ such that $\eps_n>0$, $\eps_n\to 0$ and
  \[k_0=\lim_{n\to +\infty}k_{\eps_n},\]
  then the equation
  \begin{equation}\label{eq:hj}
    -\nabla u(y)\cdot (A(y)\nabla u(y))-B(y)\cdot\nabla u(y)-k^{y}+k_0=0
  \end{equation}
%  \begin{equation}\label{eq:hj}
%    -A\abs{\nabla u(y)}^2-B\cdot\nabla u(y)-k^{y}+k_0=0
%  \end{equation}
  has a viscosity solution $u$ on $\mcy'$, belonging to $\mcc^{0,\beta}(\mcy')$ for all $\beta\in(0,1)$.
\end{lem}

\begin{proof}

%\paragraph{Step 1. Extraction of a sequence $(k_{\eps_n})_n$ which converges.}
  %First, by Lemma~\ref{lem:estimate_k}, the family $(k_{\eps})_{\eps>0}$ is bounded. It is thus possible to extract a sequence $(\eps_n)_n$ such that $k_{\eps_n}$ converges to a limit $k_{\infty}$.
  The proof is an adaptation of the proof of~\cite[Proposition 3.1]{HNR11}.
  We set:
\[v_n(y,z):=\eps_n\ln\pth{\varphi_{\eps_n}(y,z)}.\]
We use the shortcuts $\varphi_n=\varphi_{\eps_n}$ and $k_n=k_{\eps_n}$.
 We choose compact subsets $\mcy'\subset\subset\mcy$ and $\mcz'\subset\subset\mcz$ and we note $\mcw'=\mcy'\times\mcz'$.
For convenience, we assume that $0\in\mcy'$ and $0\in\mcz'$.  We renormalise $\varphi_\eps$ so that $\varphi_{\eps}(0,0)=1.$

\paragraph{Step 1. On $\overline{\mcw'}$, uniform convergence of the sequence $(v_n)_n$ to a function $v$.}
We show that up to extraction, the sequence of restricted functions $({v_n}_{|\overline{\mcw'}})_n$ converges uniformly to a function $v$ defined on $\overline{\mcw'}$. %see~\cite{HNR11}.
Consider the change of variables $h_n(y,z):=(\eps_ny,z)$ and set $\widetilde{\varphi}_n:=\varphi_n\circ h_n$.
Then $\widetilde{\varphi}_n$ is a solution of 
\[\nabla_y\cdot((A\circ h_n)\nabla_y\widetilde{\varphi}_n)+\nabla_y\cdot(\widetilde{\varphi}_n(B\circ h_n))+(\mcl_z\circ h_n+c\circ h_n)\varphi_n =0.\]
Note that the $\widetilde{\varphi}_n$ are positive, and that the coefficients are bounded uniformly in $n$.
Moreover, the diffusion coefficients $A\circ h_n$ and $a\circ h_n$ are elliptic uniformly in $n$.
By the Harnack inequality
(which, if needed, we can apply up to the boundary due to the Neumann boundary condition and the regularity of $\dr\mcw$),
there exists a constant $C_1>1$, independent of $n$, such that for each $n$, for each $y_0\in\mcy$, for each $(y,z)\in\mcw$ such that $\dabs{y-y_0}\loq\eps_n$,
\begin{equation}\label{eq:harnack_phi}
  \frac{1}{C_1}\widetilde{\varphi}_n\pth{\frac{y_0}{\eps_n},0}\loq \widetilde{\varphi}_n\pth{\frac{y}{\eps_n},z}\loq C_1\widetilde{\varphi}_n\pth{\frac{y_0}{\eps_n},0}.
\end{equation}
Moreover, by an interior estimate on the derivative (\cite[Theorem~8.32]{GilTru01}), there exists a constant $C_2>0$ (independent of $n$) such that for each $n$, for each $(y_0,z_0)\in\mcw'$,
\begin{equation}\label{eq:estimationeriv}
  \abs{\nabla\widetilde{\varphi}_n\pth{\frac{y_0}{\eps_n},z_0}}\loq C_2\sup_{(y,z)\in\mcw\vg \dabs{y-y_0}\loq \eps_n}\widetilde{\varphi}_n\pth{\frac{y}{\eps_n},z}.
\end{equation}
Using~\eqref{eq:harnack_phi} and~\eqref{eq:estimationeriv}, we conclude that for all $(y_0,z_0)\in\mcw'$,
\[\abs{\nabla v_n(y_0,z_0)}=\abs{\frac{\nabla\widetilde{\varphi}_n\pth{\frac{y_0}{\eps_n},z_0}}{\widetilde{\varphi}_n\pth{\frac{y_0}{\eps_n},z_0}}}
\loq \frac{C_2\times C_1\widetilde{\varphi}_n\pth{\frac{y_0}{\eps_n},0}}{\frac{1}{C_1}\widetilde{\varphi}_n\pth{\frac{y_0}{\eps_n},0}}
= C_1^2C_2.\]
%\loq C_1C_2.\]
This estimation being true for each $(y_0,z_0)\in\overline{\mcw'}$, we conclude that the sequence $(\nabla v_n)_n$ is bounded in $L^{\infty}(\overline{\mcw'})$.
In particular the sequence $(v_n)_n$ is also bounded in $\mcc^0(\overline{\mcw'})$.
Thanks to these facts, we infer that the sequence $(v_n)$ is bounded in $W^{1,\infty}(\overline{\mcw'})$.
%By Ascoli's theorem, the sequence $({v_n}_{|\overline{\mcw'}})_n$ is compact in $\mcc^0(\overline{\mcw'})$.
%Up to extraction, the sequence $({v_n}_{|\overline{\mcw'}})_n$ converges uniformly on $\overline{\mcw'}$ to a function $v$ defined on $\overline{\mcw'}$.
By Ascoli's theorem, for all $\beta\in(0,1)$, up to extraction, the sequence $({v_n}_{|\overline{\mcw'}})_n$ converges in $\mcc^{0,\beta}(\overline{\mcw'})$ to a function $v$ defined on $\overline{\mcw'}$.
By a diagonal extraction argument, we get a function $v$ that belongs to $\mcc^{0,\beta}(\overline{\mcw'})$ for all $\beta\in(0,1)$.

\medskip

Now we let 
\begin{align*}
  u(y):=v(y,0),\esp
  u_n(y):=v_n(y,0).
\end{align*}
Then, $u_n$ converges uniformly to $u$ on $\mcy'$, and $u$ belongs to~$\mcc^{0,\beta}(\mcy')$ for all $\beta\in(0,1)$.

\paragraph{Step 2. Find an equation satisfied by $u_n$.}
We have:
\[\nabla_yv_n=\eps_n\frac{\nabla_y\varphi_n}{\varphi_n},\]
so:
\[\eps_n\nabla_y\cdot(A\nabla_yv_n)=\eps_n^2\frac{\nabla_y\cdot(A\nabla_y\varphi_n)}{\varphi_n}-\nabla_yv_n\cdot(A\nabla_yv_n).\]
Now, by the definition of $\varphi_n$,
\begin{align*}
  \eps_n^2\frac{\nabla_y\cdot(A\nabla_y\varphi_n)}{\varphi_n}&=-\eps_n\frac{\nabla_y\cdot(\varphi_nB)}{\varphi_n}-\frac{(\mcl_z+c)\varphi_n}{\varphi_n}+k_n\\
  &=-B\cdot\nabla_yv_n-\eps_n\nabla_y\cdot B-\frac{(\mcl_z+c)\varphi_n}{\varphi_n}+k_n.
\end{align*}
We conclude that
\begin{equation*}
  -\eps_n\nabla_y\cdot(A\nabla_y v_n)-\nabla_yv_n\cdot (A\nabla_y v_n)
  -B\cdot\nabla_yv_n -\eps_n\nabla_y\cdot B
    -\frac{(\mcl_z+c)\varphi_n}{\varphi_n}+k_{n}=0.
\end{equation*}
Recalling that $u_n(y)=v_n(y,0)$, we obtain:
\begin{equation}\label{eq:visco_finite_L}
  -\eps_n\nabla\cdot(A\nabla u_n)-\nabla u_n\cdot(A\nabla u_n)
  -B\cdot\nabla u_n -\eps_n\nabla\cdot B
    -\frac{(\mcl_z+c)\varphi_n(y,0)}{\varphi_n(y,0)}+k_{n}=0.
\end{equation}

\paragraph{Step 3. Conclusion.}
Since $u$ is continuous, we have: % \cite[Chapter 6]{CIL92}}:
\[u=\pth{\limsup_{n\to+\infty}}^*u_n,\]
where, as in~\cite{CIL92},
  \[\pth{\limsup_{n\to+\infty}}^*u_n(y):=\limsup_{j\to+\infty}\acc{u_n(y')\ /\ n\goq j,\, y'\in\mcy,\, \dabs{y-y'}\loq 1/j}.\]
  Let us also define the notation
  \begin{multline*}
    J^{2,+}u(y):=\left\{(\nabla\varphi(y),\nabla^2\varphi(y))\in\Er^n\times\mcm_n(\Er)\ / \right.\\
    \left.\varphi\text{ is $\mcc^2$ near $y$ and $u-\varphi$ has a local maximum at $y$}\right\}.
  \end{multline*}
Let $y\in\mcy'$ and $(p,X)\in J^{2,+}u(y)$.
By \cite[Lemma 6.1]{CIL92},
there exist sequences $(y_n)_n$ and $(p_n,X_n)\in J^{2,+}u_n(y_n)$ such that, up to extraction of a subsequence,
\[(y_n,p_n,X_n)\to (y,p,X).\]
Note that the sequence $(X_n)_n$ is bounded and recall that $\eps_n\to 0$, so $\eps_n\Tr(A^TX_n)\to 0$.
Further, $\nabla\cdot B$ is globally bounded so $\eps_n\nabla\cdot B(y_n)\to 0$.
Moreover, Lemma~\ref{lem:cv_to_ky} yields:
\[\lim_{n\to+\infty}\frac{(\mcl_z+c)\varphi_n(y_n,0)}{\varphi_n(y_n,0)}=k^{y}.\]
Define $A_d\in\mcc^{0,\alpha}(\mcy,\Er^{n})$ by $[A_d]_i=\sum_{j=1}^n\dr_i(A_{ij})$ where $[A_d]_{i}$ is the $i^{th}$ coordinate of the vector $A_d$ and $A_{ij}$ is the coefficient with coordinates $i,j$ of the matrix $A$. Since $A_d$ is globally bounded, we have $\eps_n A_d\cdot p_n\to 0$.

We have, using the remarks above,
\begin{multline*}
  -p\cdot(A(y)p)-B(y)\cdot p-k^{y}+k_0\\
  =\lim_{n\to+\infty}\Big[\pth{-\eps_n\pth{\Tr(A^TX_n)+A_d\cdot p_n}-p_n\cdot(Ap_n)-B\cdot p_n-\eps_n\nabla\cdot B}(y_n)\\-\pth{\frac{(\mcl_{z}+c)\varphi_n}{\varphi_n}}(y_n,0)+k_n\Big].
\end{multline*}
Now, the first line of ~\eqref{eq:visco_finite_L} and the fact that $(p_n,X_n)\in J^{2,+}u(y_n)$ imply: for all $n\goq 0$,
\begin{multline*}
  \pth{-\eps_n\pth{\Tr(A^TX_n)+A_d\cdot p_n}-p_n\cdot (Ap_n)-B\cdot p_n-\eps_n\nabla\cdot B}(y_n)\\-\pth{\frac{(\mcl_{z}+c)\varphi_n}{\varphi_n}}(y_n,0)+k_n\loq 0.
\end{multline*}
%\begin{align*}
%  0&=-\eps_n\nabla_y\cdot\pth{A\nabla_yv_n}-A\abs{\nabla_yv_n}^2-B\cdot\nabla_yv_n-c-\frac{\nabla_{z}\cdot\pth{\mu\nabla_{z}\varphi_n}}{\varphi_n}+k_{n}\\
%  &=-\eps_n\pth{A\Delta_yv_n+\nabla_yA\cdot\nabla_yv_n}-A\abs{\nabla_yv_n}^2-B\cdot\nabla_yv_n-c-\frac{\nabla_{z}\cdot\pth{\mu\nabla_{z}\varphi_n}}{\varphi_n}+k_{n}.
%\end{align*}
We obtain, therefore:
\begin{align*}
  -p\cdot (Ap)-B\cdot p-k^{y}+k_0&\loq 0.
\end{align*}
Thus $u$ is a viscosity subsolution of~\eqref{eq:hj}. Likewise, defining $\liminf_*$ analogously to $\limsup^*$,
\[u(y)=\pth{{\liminf_{n\to+\infty}}}_{*}u_n(y).\]
Thus $u$ is a viscosity supersolution of~\eqref{eq:hj}. Therefore, $u$ is a viscosity solution of~\eqref{eq:hj}.
\end{proof}

  \begin{req}\label{req:conjecture_phi_slow}
    %We are now ready to understand the remarks on the behaviour of
    %$\frac{\varphi_{\eps}(\cdot,z_0)}{\varphi_{\eps}(\mcy,z_0)}$
    %as $\eps\to 0$ that we mentioned in the introduction.
    The proof of Lemma~\ref{lem:visco_sol} suggests that as $\eps\to0$,
    \[\varphi_{\eps}(y,0)\simeq e^{u(y)/\eps},\]
    where $u$ solves~\eqref{eq:hj}.
    Therefore, we expect that
    \[\frac{\varphi_{\eps}(y,0)}{\varphi_{\eps}(y,\mcy)}\to0\]
    everywhere but at points $y_0\in\overline{\mcy}$ such that $u$ is maximal at $y_0$.
    If $u$ is of class $\mcc^1$, then the equation
    \[-\nabla u(y)\cdot (A(y)\nabla u(y))-B(y)\cdot\nabla u(y)-k^{y}+k_0=0\]
    holds in the classical sense.
    If moreover $\mcy$ has empty boundary, then all maximum points $y_0$ of $u$ are interior, so we have $\nabla u(y_0)=0$ and
    therefore $k^{y_0}=k_0$.
    Thus the limit of
    $\frac{\varphi_{\eps}(y,0)}{\varphi_{\eps}(y,\mcy)}$
    should be supported in
    $\Gamma:=\acc{y\in\mcy\ /\ k^{y}=k_0}$.
    Since $\frac{\varphi_{\eps}(y,0)}{\varphi_{\eps}(y,\mcy)}$
    has total mass $1$, it is natural to expect that 
    there exists a measure~$\mu$ on~$\Gamma$,
    with total mass equal to $1$,
    such that for all $f\in\mcc^{\infty}(\mcy)$,
    \[\int_{\mcy}f(y)\frac{\varphi_{\eps}(y,z_0)}{\varphi_{\eps}(\mcy,z_0)}\de y\to\int_{\Gamma} f(y)\mu(\de y).\]
    This holds also if we replace~$0$ with any $z_0\in\mcz$.
  \end{req}

\subsection{Main proofs}

We first prove Theorem~\ref{thm:bzero} (where it is assumed that $B\equiv0$) by studying the Hamilton-Jacobi equation given by Lemma~\ref{lem:visco_sol}.
Next, we prove Theorem~\ref{thm:general} by solving explicitly the Hamilton-Jacobi equation in the periodic $1$-dimensional setting.

\begin{proof}[Proof of Theorem~\ref{thm:bzero}]% in the general case]
  Assume that $B\equiv 0$. Take $\eps_n\to 0$ such that $k_{\eps_n}$ converges to some $k_0\in\Er$. Lemma~\ref{lem:visco_sol} then implies that for all $\mcy'\subset\subset\mcy$, there exists a viscosity solution $u\in\mcc^{0,\beta}(\mcy')$ (for all $\beta\in(0,1)$)
  of the equation
  \[-\nabla u\cdot (A\nabla u)-k^{y}+k_0=0.\]
  Starting from this point, the proof is cut into two steps.
%  This implies in particular that for all $y\in\mcy$, %rq: {no pb: $J^{2,+}$ empty on a nb implies the function is monotonous!}
  %  $k_{0}\goq k^{y}.$

  \paragraph{Step 1. For all $y\in\mcy$, $k_0\goq k^y$.}
  Let $y\in\mcy$ and let us show that $k_{0}\goq k^{y}$.
    First, assume that there exists a function $\phi$ which is $\mcc^2$ in a neighbourhood of $y$ and such that $u-\phi$ has a local minimum at $y$.
    Then, we have $-\nabla\phi\cdot(A\nabla\phi)-k^{y}+k_0\goq 0$ and, owing to the ellipticity of $A$,
    we have $\nabla\phi\cdot(A\nabla\phi)\goq 0$.
    Thus $k_{0}\goq k^{y}$.
    
    If there is no such $\phi$,
    we construct a sequence $y_n\to y$ such that $y_n$ is associated with a $\phi_n$,
    and apply the same reasoning at $y_n$.
    Namely, we use the $\mcc^{0,1/2}$-regularity of $u$ to define:
    \[a_n:=\sup_{y'\in\mcy,\,0<\dabs{y'-y}<1/n}\frac{\abs{u(y')-u(y)}}{\sqrt{\dabs{y'-y}}}.\]
%    \[a_n:=\inf\acc{a>0 \ /\ \text{for all } y'\in \mcy\text{ with } \dabs{y'-y}<\frac{1}{n},\ u(y')\goq u(y)-a\sqrt{\dabs{y-y'}}}\in(0,+\infty).\]
    We set:
    \[b_n:=\inf_{y'\in\mcy,\,\dabs{y'-y}<1/n}\pth{u(y')+4{a_nn^{3/2}}\dabs{y'-y}^2}\]
%    \[b_n:=\sup\acc{b\in\Er \ /\ \text{for all } y'\in \mcy\text{ with } \dabs{y'-y}<\frac{1}{n},\ u(y')\goq b+u(y)-a_n\dabs{y-y'}^2}\in(0,+\infty).\]
    and we let
    \[\phi_n(y'):=b_n-4{a_nn^{3/2}}\dabs{y'-y}^2.\]
    We note that $u(y)\goq b_n$ and that $\phi_n(y')\loq u(y')$ for $\dabs{y-y'}\loq 1/n$.
    Moreover, there exists~$y_n$ with $\dabs{y_n-y}\loq 1/n$ such that
    $\phi_n(y_n)= u(y_n)$.
    Last, for $1/(2n)<\dabs{y'-y}\loq 1/n$,
    \begin{align*}
      u(y')&\goq u(y)-a_n\sqrt{\dabs{y'-y}}
      \goq u(y)-\frac{a_n}{\sqrt{n}},
    \end{align*}
    so, using $1/(2n)<\dabs{y'-y}$,
    \begin{align*}
      u(y')
      &> u(y)-4{a_nn^{3/2}}\dabs{y'-y}^2\\
      &\goq b_n-4{a_nn^{3/2}}\dabs{y'-y}^2=\phi_n(y').
    \end{align*}
    To summarise, we have: $\dabs{y_n-y}\loq 1/n$, $\phi_n(y_n)=u(y_n)$ and
    \begin{align*}
      %\phi_n(y_n)&=u(y_n),\\
      \phi_n(y')&\loq u(y')&\text{for all $\dabs{y'-y}\loq 1/n$},\\
      \phi_n(y')&< u(y')&\text{for all $1/(2n)<\dabs{y'-y}\loq 1/n$}.
    \end{align*}
    %In particular, $\dabs{y'-y}\loq 2/n$, so the inequality $\phi_n\loq u$
    %holds in a neighbourhood of $y_n$.
    %Thus, for all $n$, the function
    Therefore, $\dabs{y_n-y}\loq 1/(2n)$ and $u-\phi_n$ has a local minimum at~$y_n$.
    Since~$\phi_n$ is~$\mcc^2$, we conclude as above that $k_{0}\goq k^{y_n}$.
    We get, using $y_n\to y$ and the continuity of $y\mapsto k^y$, that $k_{0}\goq k^{y}$.

    \paragraph{Step 2. There exists $y_{\infty}\in\overline{\mcy}$ such that $k_0=k^{y_{\infty}}$.}

    Recall that the functions $u_n=\eps_n\ln(\varphi_{\eps_n}(\cdot,0))$ have been defined in the proof of Lemma~\ref{lem:visco_sol}. For $n\goq 0$, let $y_n\in\overline{\mcy}$ satisfy
  \[u_n(y_n)=\max_{y\in\overline{\mcy}}u_n(y).\]
  Since $\overline{\mcy}$ is compact, the sequence $(y_n)_n$ converges, up to extraction, to some $y_{\infty}\in\overline{\mcy}$. By Lemma~\ref{lem:cv_to_ky}, maybe up to another extraction, we have
  \[\frac{(\mcl_z+c)\varphi_n(y_n,0)}{\varphi_n(y_n,0)}\to k^{y_{\infty}}.\]
  If $y_n\in\mcy$, then $\nabla u_n(y_n)=0$. If $y_n\in\dr\mcy$, then the condition $\nu\cdot A\nabla u_n=0$ % extrema liés
  on $\dr\mcy$ implies that $\nabla u_n(y_n)=0$. In either case,~\eqref{eq:visco_finite_L} implies that for all $n\goq 0$,
  \[-\eps_n\nabla\cdot(A\nabla u_n)(y_n)-\eps_n\nabla\cdot B
  -\frac{(\mcl_z+c)\varphi_n(y_n,0)}{\varphi_n(y_n,0)}+k_{n}=0.
  \]
  Taking $n\to+\infty$ gives (using the same estimates as in Step~3 of the proof of Lemma~\ref{lem:visco_sol}): $-k^{y_{\infty}}+k_0=0$. Hence $k_0=k^{y_{\infty}}$ so necessarily
  \[k_0=\max_{y\in\overline{\mcy}}k^y.\]
  Finally, the limit is unique, so the whole family $(k_{\eps})_{\eps>0}$ converges to $k_0$.
\end{proof}

Now, we prove Theorem~\ref{thm:general}: we drop the assumption that $B\equiv 0$, but we assume that $\mcy$ is the $1$-dimensional torus.
There are two cases: either $\abs{\gamma}\goq j(M)$, either $\abs{\gamma}< j(M)$.
We first prove the former.

\begin{proof}[Proof of Theorem~\ref{thm:general}, case $\abs{\gamma}\goq j(M)$.]
  The proof is directly inspired from Section 3 of~\cite{HNR11}. 
  We assume that $\mcy$ is the $1$-dimensional torus. Then $\mcy\subset\subset\mcy$, so we can take $\mcy'=\mcy$ in Lemma~\ref{lem:visco_sol}.
We first prove that there exists only one $k_0$ such that ~\eqref{eq:hj} has a solution on $\mcy$, \ie a $1$-periodic viscosity solution defined on~$\Er$. Then, we make explicit such a solution and conclude.

  \paragraph{Step 1. Uniqueness of $k_0$ such that~\eqref{eq:hj} has a periodic solution.}
  We assume that there exist $k_{0}\in\Er$ and $\hat{k}_{0}>k_{0}$ such that there are periodic viscosity solutions $u$ and $\hat{u}$ to:
\begin{align*}
    -A(y)(u'(y))^2-B(y) u'(y)-k^{y}+k_0&=0,\\
    -A(y)(\hat{u}'(y))^2-B(y)\hat{u}'(y)-k^{y}+\hat{k}_0&=0.
\end{align*}
%\begin{align*}
%    -\nabla u(y)\cdot(A(y)\nabla u(y))-B(y)\cdot\nabla u(y)-k^{y}+k_0&=0,\\
%    -\nabla \hat{u}(y)\cdot(A(y)\nabla \hat{u}(y))-B(y)\cdot\nabla \hat{u}(y)-k^{y}+\hat{k}_0&=0.
%\end{align*}
Since $k_{0}\neq \hat{k}_{0}$, we must have $u\neq\hat{u}$. 
Therefore, up to adding a constant to $u$, we may assume that there exists an interval $U\subset\Er$ such that:
\begin{equation}\label{eq:condition_u}
  \left\{
  \begin{aligned}
    u(y)&>\hat{u}(y)&\text{at some point $y\in U$,}\\
    u&=\hat{u}&\text{on $\dr U$}.
  \end{aligned}
  \right.
\end{equation}
But, since $\hat{k}_{0}>k_{0}$, there is a small $\eps>0$ such that on $U$:
\begin{align*}
  &\eps u+A(y)(u'(y))^2+B(y) u'(y)\\
  &\esp\loq\eps \hat{u}+ A(y)(\hat{u}'(y))^2+B(y)\hat{u}'(y)
\end{align*}
%\begin{align*}
%  &\eps u+\nabla u(y)\cdot(A(y)\nabla u(y))+B(y)\cdot\nabla u(y)\\
%  &\esp\loq\eps \hat{u}+\nabla u(y)\cdot(A(y)\nabla u(y))+B(y)\cdot\nabla \hat{u}(y)
%\end{align*}
in the viscosity sense. 
Now, with the second line of~\eqref{eq:condition_u}, we can apply the maximum principle, which implies that $u\loq\hat{u}$ in $U$. We reach a contradiction with the first line of~\eqref{eq:condition_u}: this ensures the uniqueness of the limit of $k_{\eps}$ as $\eps\to 0$.

By Lemma~\ref{lem:estimate_k}, the family $(k_{\eps})_{\eps>0}$ is bounded. Hence the whole family $(k_\eps)_\eps$ converges to~$k_0$ as $\eps\to 0$:
\[k_0=\lim_{\eps\to 0}k_{\eps}.\]
%\paragraph{Step 2. Lower bound on $k_0$.} 
Let us rewrite~\eqref{eq:hj} in the following way:
%\[-A(y)\pth{\nabla u(y)-\frac{B(y)}{2A(y)}}^2-\pth{k^{y}-\frac{{B(y)}^2}{4A(y)}}+k_0=0.\]
\begin{equation}\label{eq:hj2}
  -A\pth{ u'+\frac{B}{2A}}^2-\pth{k^{y}-\frac{{B}^2}{4A}}+k_0=0.
\end{equation}
%This implies in particular that for all $y\in\mcy$, 
%\[k_0\goq {k^{y}-\frac{{B}(y)^2}{4A(y)}},\]
%so we must have $k_0\goq M=\max\pth{k^y-\frac{B^2(y)}{4A(y)}}$. The first statement of the corollary is proved.
\paragraph{Step 2. Conclusion.}
Recall that
\[\gamma=\int_0^1\frac{B(y)}{2A(y)}\de y,\esp j(M)=\int_0^1\sqrt{\frac{M-k^y+\frac{B(y)^2}{4A(y)}}{A(y)}}\de y.\]
If $\gamma\goq j(M)$, then there exists $k\goq M$ such that
  \begin{equation}\label{eq:k_wrt_lambda}
    \int_0^1\frac{B(y)}{2A(y)}\de y=j(k)=\int_0^1\sqrt{\frac{k-\pth{k^{y}-\frac{{B(y)^2}}{4A(y)}}}{A(y)}}\de y.
  \end{equation}
  Then, the function 
%  \[u(x):=\frac{Bx}{2A}-\int_0^x\sqrt{\frac{k-\pth{k^{y}-\frac{{B^2}}{4A}}}{A}}\de y\]
  \[\overline{u}(y):=-\int_0^y\frac{B(y')}{2A(y')}\de y'+
  \int_0^y\sqrt{\frac{k-\pth{k^{y'}-\frac{{B(y')^2}}{4A(y')}}}{A(y')}}\de y'\]
  is $1$-periodic. Moreover, the function $\overline{u}$ solves~\eqref{eq:hj2} (with $k$ instead of $k_0$) and, therefore, solves~\eqref{eq:hj}. By Lemma~\ref{lem:visco_sol} and the uniqueness stated in Step 1, we conclude that $k_{0}$ and $k$ must coincide, so that $k_{0}$ satisfies~\eqref{eq:k_wrt_lambda}:
  \[\int_0^1\frac{B(y)}{2A(y)}\de y=j(k_{0}).\]
  This concludes the proof for $\gamma\goq j(M)$. If $\gamma\loq -j(M)$, then there exists $k\goq M$ such that
  \begin{equation*}
    -\int_0^1\frac{B(y)}{2A(y)}\de y=j(k)=\int_0^1\sqrt{\frac{k-\pth{k^{y}-\frac{{B(y)^2}}{4A(y)}}}{A(y)}}\de y.
  \end{equation*}
  We conclude in the same way, with the function
  \[\overline{u}_{\gamma<0}(y):=\int_0^y\frac{B(y')}{2A(y')}\de y'+
  \int_0^y\sqrt{\frac{k-\pth{k^{y'}-\frac{{B(y')^2}}{4A(y')}}}{A(y')}}\de y'\]
  instead of $\overline{u}$.
\end{proof}

\begin{proof}[Proof of Theorem~\ref{thm:general}, case $\abs{\gamma}\loq j(M)$.]
  First, we consider $0\loq\gamma<j(M)$, \ie
  \[0\loq\int_0^1\frac{B(y)}{2A(y)}\de y< \int_0^1\sqrt{\frac{M-k^y+\frac{B(y)^2}{4A(y)}}{A(y)}}\de y.\]
  We assume for convenience that $M$ is reached at $0$, \ie
  \[M=\max_{y\in\mcy}\pth{k^y-\frac{B(y)^2}{4A(y)}}=k^0-\frac{B(0)^2}{4A(0)}.\]
  Since all these functions are $1$-periodic, $M$ is also reached at each $y\in\Zed$.

  \paragraph{Step 1. Construct a viscosity solution of an associated problem.}

  We define a $2$-periodic function $S$ by
  \[S(y):=\left\{
  \begin{aligned}
    1,&&y\in\ifo{0,1},\\
    -1,&&y\in\ifo{1,2}.
  \end{aligned}
  \right.\]
  Now, we define $v$ by:
  \[v(y):=-\int_0^y\frac{B(y')}{2A(y')}\de y'+\int_0^yS(y')\sqrt{\frac{M-\pth{k^{y'}-\frac{{B(y')^2}}{4A(y')}}}{A(y')}}\de y'.\]
%  Note that $v(2)=-2\gamma$, so 
%  $y\mapsto v(y)+\gamma y$
%  is $2$-periodic.
  For $y_0\in\Er\setminus\Zed$, the function $S$ is continuous at $y_0$ so $v$ is differentiable at $y_0$,  and satisfies, at $y=y_0$,
  \begin{equation}\label{eq:hj_M}
    -A(v')^2-Bv'-k^{y}+M=0.
  \end{equation}
  For $y_0\in\Zed$, we have
  \[M-\pth{k^{y_0}-\frac{B(y_0)^2}{4A(y_0)}}=0,\]
  so $v$ is also differentiable at $y_0$ and satisfies~\eqref{eq:hj_M} at $y=y_0$ in the classical sense.
  In particular, $v$ is a viscosity solution of~\eqref{eq:hj_M}.

  \paragraph{Step 2. Prove a uniqueness result and conclude.}
  Let $k_0\in\Er$ be a limit of a subsequence of $(k_{\eps})_{\eps}$ and let~$u$ be given by Lemma~\ref{lem:visco_sol}.
  We have $v(0)=0$ and $v(1)=j(M)-\gamma>0$.
  We add a constant to $u$ in such a way that $u(0)\in (0, v(1))$; Equation~\eqref{eq:hj} remains satisfied by $u$.
  
  We have, by $1$-periodicity of $u$:
  \[u(0)>v(0),\qquad u(1)<v(1).\]
  We note that
  \[y\mapsto v(y)+\gamma y\]
  is $2$-periodic; therefore, we have $v(y)\to+\infty$ as $y\to-\infty$ and $v(y)\to-\infty$ as $y\to+\infty$.
  Thanks to these elements, we get the existence of $y_0<0$, $y_1\in(0,1)$ and $y_2>1$ such that $u(y_i)=v(y_i)$ for all $i=1,2,3$.
  Therefore, there exist two bounded intervals $U:=(y_0,y_1)$ and $U':=(y_1,y_2)$ such that:
  \begin{equation}\label{eq:intervals_uup}
  \left\{
  \begin{aligned}
    u(y)&>v(y)&\text{at some point $y\in U$,}\\
    v(y)&>u(y)&\text{at some point $y\in U'$,}\\
    u&=v&\text{on $\dr U$ and $\dr U'$}.
  \end{aligned}
  \right.
  \end{equation}
  Recall that $v$ is a viscosity solution of~\eqref{eq:hj_M}
  and that $u$ is a viscosity solution of~\eqref{eq:hj_M} but with $k_0$ instead of $M$.
  Using the same reasoning as in Step~1 of the previous proof (starting from~\eqref{eq:condition_u}),
  the first and last line of~\eqref{eq:intervals_uup} imply that $M\loq k_0$;
  the second and last line of~\eqref{eq:intervals_uup} imply that $M\goq k_0$.
  Hence, $k_0=M$.
  
  \medskip
  
  Finally, the case $-j(M)<\gamma\loq 0$ can be proved in the same way, by considering
  \[v_{\gamma<0}(y):=\int_0^y\frac{B(y')}{2A(y')}\de y'+\int_0^yS(y')\sqrt{\frac{M-\pth{k^{y'}-\frac{{B(y')^2}}{4A(y')}}}{A(y')}}\de y'\]
  instead of $v$.
\end{proof}

\section*{Acknowledgements}
The author thanks Raphaël Forien, François Hamel and Lionel Roques for their advice and support.
The author also thanks the reviewers whose helpful comments allowed him to improve the manuscript, and Nicolas Champagnat, who pointed out a mistake in the first version of the proof of Proposition~\ref{ppn:estimation_diff_zero}. \\
This work was supported by the French Agence Nationale de la Recherche (ANR-18-CE45-0019 `RESISTE' and {ANR-23-CE40-0023-01} 'ReaCH') and by the Chaire Modélisation Mathématique et Biodiversité (École Polytechnique, Muséum national d’Histoire naturelle, Fondation de l’École Polytechnique, VEOLIA Environnement).

\printbibliography

\end{document}